\input amstex
\documentstyle{amsppt}
\loadbold
\magnification=\magstep1
\pageheight{9.0truein}
\pagewidth{6.5truein}
\NoBlackBoxes
\TagsOnLeft

\input xy
\xyoption{matrix}\xyoption{arrow}

\def\pup#1{{\rm(}#1{\rm)}}

\def\ZZ{{\Bbb Z}}

\def\NN{{\Bbb N}}
\def\CC{{\Bbb C}}

\def\O{{\Cal O}}

\def\bfp{{\bold p}}
\def\bfq{{\bold q}}
\def\bfla{{\boldsymbol\lambda}}
\def\bfmu{{\boldsymbol\mu}}
\def\bfxi{{\boldsymbol\xi}}
\def\bfsig{{\boldsymbol\sigma}}
\def\bfz{{\bold 0}}

\def\Ahat{\widehat{A}}

\def\alphahat{\widehat{\alpha}}

\def\gfrak{{\frak g}}
\def\hfrak{{\frak h}}
\def\mfrak{{\frak m}}
\def\glfrak{{\frak {gl}}}

\def\hsl{\hbar}

\def\eps{\epsilon}

\def\Lie{\operatorname{Lie}}

\def\Fract{\operatorname{Fract}}

\def\spec{\operatorname{spec}}

\def\Pspec{\operatorname{P{.}spec}}

\def\Der{\operatorname{Der}}
\def\Aut{\operatorname{Aut}}
\def\End{\operatorname{End}}
\def\kx{k^\times}

\def\OlpMn{\O_{\lambda,\bfp}(M_n(k))}
\def\khh{k[[\hsl]]}
\def\OMn{\O(M_n(k))}
\def\lexle{<_{\operatorname{lex}}}
\def\sign{\operatorname{sign}}
\def\tr{\operatorname{tr}}
\def\chr{\operatorname{char}}
\def\adj{\operatorname{adj}}
\def\trdeg{\operatorname{tr{.}deg}}

\def\AlDu{{\bf 1}}
\def\BrGo{{\bf 2}}
\def\BrGr{{\bf 3}}
\def\Cau{{\bf 4}}
\def\Cohn{{\bf 5}}
\def\DiLe{{\bf 6}}
\def\pdixmo{{\bf 7}}
\def\GLet{{\bf 8}}
\def\specstrat{{\bf 9}}
\def\GoYa{{\bf 10}}
\def\Hrt{{\bf 11}}
\def\Kam{{\bf 12}}
\def\KoWa{{\bf 13}}
\def\Loo{{\bf 14}}
\def\Mus{{\bf 15}}
\def\New{{\bf 16}}
\def\Nou{{\bf 17}}
\def\Ohcat{{\bf 18}}
\def\Oh{{\bf 19}}
\def\Ohtwo{{\bf 20}}
\def\Ohthree{{\bf 21}}
\def\Pan{{\bf 22}}
\def\Ric{{\bf 23}}
\def\St{{\bf 24}}
\def\TauYu{{\bf 25}}
\def\TaYu{{\bf 26}}
\def\Van{{\bf 27}}
\def\Ver{{\bf 28}}

\topmatter

\title The Dixmier-Moeglin Equivalence and a Gel'fand-Kirillov Problem for Poisson
Polynomial Algebras
\endtitle

\rightheadtext{Poisson polynomial rings}

\author K. R. Goodearl and S. Launois \endauthor

\address Department of Mathematics, University of California, Santa Barbara, CA 93106, USA
\endaddress

\email goodearl\@math.ucsb.edu \endemail

\address University of Edinburgh and Maxwell Institute for Mathematical Sciences,
School of Mathematics, JCMB, King's Buildings, Mayfield Road, Edinburgh EH9 3JZ, Scotland
\endaddress

\curraddr (From 01 August 2007) Institute of Mathematics, Statistics and Actuarial Science,
University of Kent at Canterbury, CT2 7NF, UK \endcurraddr

\email stephane.launois\@ed.ac.uk \endemail

\thanks This research was partially supported by Leverhulme Research Interchange Grant
F/00158/X (UK), that of the first author by a grant from the National Science Foundation (USA),
and that of the second author by a Marie Curie Intra-European
 Fellowship within the $6^{\text{th}}$ European Community Framework
Programme.
\endthanks

 \subjclass Primary 17B63 \endsubjclass

\abstract The structure of Poisson polynomial algebras of the type obtained as semiclassical
limits of quantized coordinate rings is investigated. Sufficient conditions for a rational Poisson
action of a torus on such an algebra to leave only finitely many Poisson prime ideals invariant
are obtained. Combined with previous work of the first-named author, this establishes the Poisson
Dixmier-Moeglin equivalence for large classes of Poisson polynomial rings, such as semiclassical
limits of quantum matrices, quantum symplectic and euclidean spaces, quantum symmetric and
antisymmetric matrices. For a similarly large class of Poisson polynomial rings, it is proved that
the quotient field of the algebra (respectively, of any Poisson prime factor ring) is a rational
function field $F(x_1,\dots,x_n)$ over the base field (respectively, over an extension field of
the base field) with $\{x_i,x_j\}= \lambda_{ij} x_ix_j$ for suitable scalars $\lambda_{ij}$, thus
establishing a quadratic Poisson version of the Gel'fand-Kirillov problem. Finally, partial
solutions to the isomorphism problem for Poisson fields of the type just mentioned are obtained.
\endabstract

\endtopmatter

\document

\head 0. Introduction \endhead

Fix a base field $k$ of characteristic zero throughout. All algebras are assumed
to be over $k$, and all relevant maps (automorphisms, derivations, etc.) are
assumed to be $k$-linear.

Recall that a {\it Poisson algebra\/} (over $k$) is a commutative $k$-algebra $A$
equipped with a Lie bracket $\{-,-\}$ which is a derivation (for the associative
multiplication) in each variable. We investigate ({\it iterated\/}) {\it Poisson
polynomial algebras\/} over $k$, that is, polynomial algebras $k[x_1,\dots,x_n]$
equipped with Poisson brackets such that
$$\bigl\{ x_i,\, k[x_1,\dots,x_{i-1}] \bigr\} \subseteq k[x_1,\dots,x_{i-1}]x_i +
k[x_1,\dots,x_{i-1}]$$
for $i=2,\dots,n$ (see \S1.1 for more detail on the conditions satisfied by
such a bracket). Many such Poisson algebras are semiclassical limits of quantum algebras, and
these provide our motivation and focus (see Section 2). The Kirillov-Kostant-Souriau Poisson
structure on the symmetric algebra of a finite dimensional Lie algebra $\gfrak$ can be put in the
form of a Poisson polynomial algebra when $\gfrak$ is completely solvable. This also holds for the
basic example of a {\it Poissson-Weyl algebra\/}, namely a polynomial algebra
$k[x_1,\dots,x_n,y_1,\dots,y_n]$ equipped with the Poisson bracket such that
$$\xalignat2 \{x_i,x_j\} &= \{y_i,y_j\} =0  &\{x_i,y_j\} &= \delta_{ij}  \tag0-1
\endxalignat$$
for all $i$, $j$.

Our investigation has two main goals, namely to establish conditions under which
Poisson analogs of the Dixmier-Moeglin equivalence and (a quadratic analog of) the
Gel'fand-Kirillov problem hold for Poisson polynomial rings.

\definition{0.1\. The Poisson Dixmier-Moeglin equivalence} Let $A$ be a Poisson
algebra. A {\it Poisson ideal\/} of $A$ is any ideal $I$ such that $\{A,I\}
\subseteq I$, and a {\it Poisson prime\/} ideal is any prime ideal which
is also a Poisson ideal. The set of Poisson prime ideals in $A$ forms the {\it
Poisson prime spectrum\/}, denoted $\Pspec A$, which is given the relative Zariski
topology inherited from $\spec A$. Given an arbitrary ideal $J$ of $A$, there is a
largest Poisson ideal contained in $J$, called the {\it Poisson core\/} of $J$.
The {\it Poisson primitive\/} ideals of $A$ are the Poisson cores of the maximal
ideals. (One thinks of the Poisson core of an ideal in a Poisson algebra as
analogous to the {\it bound\/} of a left ideal $L$ in a noncommutative algebra
$B$, that is, the largest two-sided ideal of $B$ contained in $L$.) The Poisson primitive ideals in the coordinate ring of a complex affine Poisson variety $V$ are the defining ideals of the Zariski closures of the symplectic leaves in $V$ \cite{\BrGo, Lemma 3.5}, and they are the key to Brown and Gordon's concept of {\it symplectic cores\/} \cite{\BrGo, \S3.3}.

The {\it Poisson center\/} of $A$ is the subalgebra
$$Z_p(A)= \{z\in A\mid \{z,-\} \equiv 0\}.$$
For any Poisson prime ideal $P$ of $A$, there is an induced Poisson bracket on
$A/P$, which extends uniquely to the quotient field $\Fract A/P$ (e.g.,
\cite{\Loo, Proposition 1.7}). We say that $P$ is {\it Poisson rational\/} if the
field $Z_p(\Fract A/P)$ is algebraic over $k$.

By analogy with the Dixmier-Moeglin equivalence for enveloping algebras, we say
that $A$ satisfies the {\it Poisson Dixmier-Moeglin equivalence\/} (e.g.,
\cite{\Oh, pp\. 7,8}) provided the following sets coincide:
\roster
\item The set of Poisson primitive ideals in $A$;
\item The set of locally closed points in $\Pspec A$;
\item the set of Poisson rational Poisson prime ideals of $A$.
\endroster
If $A$ is an affine (i.e., finitely generated) $k$-algebra, then $(2)\subseteq
(1)\subseteq (3)$ \cite{\Oh, Propositions 1.7, 1.10}, so the main difficulty is
whether $(3)\subseteq (2)$. No examples are known of affine Poisson algebras for
which the Poisson Dixmier-Moeglin equivalence fails. The equivalence has been
established in \cite{\pdixmo} for Poisson algebras with suitable torus actions, as
follows.
\enddefinition

\definition{0.2\. Torus actions} Suppose that $H$ is a group acting on a Poisson
algebra $A$ by {\it Poisson automorphisms\/} (i.e., $k$-algebra automorphisms that
preserve the Poisson bracket). For each $H$-stable Poisson prime $J$ of $A$, set
$$\Pspec_J A= \{P\in \Pspec A\mid \bigcap_{h\in H} h(P)=J\},$$
the {\it $H$-stratum\/} of $\Pspec A$ corresponding to $J$. These $H$-strata
partition $\Pspec A$ as $J$ runs through the $H$-stable Poisson primes of $A$.

Now assume that $H= (\kx)^r$ is an algebraic torus over $k$. In this case, the action of $H$ on
$A$ is called {\it rational\/} provided $A$ is generated (as a $k$-algebra) by $H$-eigenvectors
whose eigenvalues are rational characters of $A$. (See \S1.4 for the general definition of a
rational action of an algebraic group, and \cite{\BrGo, Theorem II.2.7} for the equivalence with
the above condition in the case of a torus.) Rationality will be clearly satisfied for the torus
actions given in the  examples in Section 2. In view of the following theorem, all we will need to
establish is that the number of $H$-stable Poisson primes is finite in these examples.
\enddefinition

\proclaim{0.3\. Theorem} {\rm\cite{\pdixmo, Theorem 4.3}} Let $A$ be an affine
Poisson algebra and $H= (\kx)^r$ an algebraic torus acting rationally on $A$ by
Poisson automorphisms. Assume that $A$ has only finitely many $H$-stable Poisson
prime ideals.

Then the Poisson Dixmier-Moeglin equivalence holds in $A$, and the Poisson
primitive ideals are precisely those Poisson primes maximal in their $H$-strata.
\qed\endproclaim

\definition{0.4\. A quadratic Poisson Gel'fand-Kirillov problem} The original
Gel'fand-Kirillov problem asked whether the quotient division ring of the enveloping algebra of a
finite dimensional algebraic Lie algebra $\gfrak$ over $k$ is isomorphic to a Weyl skew field over
a purely transcendental extension $K$ of $k$, i.e., the quotient division ring of a Weyl algebra
over $K$. Vergne raised the corresponding question for the Kirillov-Kostant-Souriau Poisson
structure on the symmetric algebra of $\gfrak$, namely whether the quotient field of $S(\gfrak)$
is isomorphic (as a Poisson algebra) to the quotient field of a Poisson-Weyl algebra \cite{\Ver,
Introduction}, and answered this positively for nilpotent $\gfrak$ \cite{\Ver, Th\'eor\`eme 4.1}.
We shall use the term {\it Poisson-Weyl field\/} for the quotient field of a Poisson-Weyl algebra,
that is, for a rational function field $K(x_1,\dots,x_n,y_1,\dots,y_n)$ equipped with the (unique)
$K$-linear Poisson bracket satisfying (0-1). Vergne's result was extended to algebraic solvable
Lie algebras $\gfrak$, and to Poisson prime factors of $S(\gfrak)$ for such $\gfrak$, by Tauvel
and Yu \cite{\TaYu, Corollaire 11.8}. We also mention that Kostant and Wallach showed that a Galois extension of the
quotient field of $\O(M_n(\CC))$, with a natural Poisson structure, is a Poisson-Weyl field
\cite{\KoWa, Theorem 5.24}.

The above form of the Poisson Gel'fand-Kirillov problem, however, is not appropriate for the
algebras we consider. In fact, as we prove in Corollary 5.3, the quotient field of a semiclassical
limit of a typical quantum algebra can never be isomorphic to a Poisson-Weyl field. A suitable
version is suggested by quantum results, as follows.

Quantum versions of the Gel'fand-Kirillov problem have been studied by a number of
authors (e.g., see \cite{\BrGo, pp\. 230-231} for a summary). These involve
quotient division rings of quantized Weyl algebras, which turn out to be
isomorphic to quotient division rings of {\it quantum affine spaces\/}
$$\O_{\bfq}(k^n)= k\langle x_1,\dots,x_n\mid x_ix_j= q_{ij}x_jx_i \text{\ for all\
} i,j\rangle$$
for multiplicatively antisymmetric matrices $\bfq= (q_{ij}) \in M_n(\kx)$.
Semiclassical limits of quantum affine spaces are Poisson polynomial rings
$k[x_1,\dots,x_n]$ with Poisson brackets satisfying
$$\{x_i,x_j\}= \lambda_{ij}x_ix_j  \tag 0-2$$
for all $i$, $j$, where $\bfla= (\lambda_{ij})$ is an antisymmetric $n\times n$
matrix over $k$ (see \S2.2). It will be convenient to denote this Poisson
polynomial algebra by $k_{\bfla}[x_1,\dots,x_n]$, the corresponding Poisson
Laurent polynomial ring by $k_{\bfla}[x_1^{\pm1},\dots,x_n^{\pm1}]$, and the
corresponding Poisson field by $k_{\bfla}(x_1,\dots,x_n)$. In all three cases, the
Poisson bracket is uniquely determined by (0-2), for instance because Poisson
brackets extend uniquely to localizations \cite{\Loo, Proposition 1.7}. In the
present situation, however, the extensions are easier to establish, since we can
give them by the formula
$$\{f,g\}= \sum_{i,j=1}^n \lambda_{ij}x_ix_j \dfrac{\partial f}{\partial
x_i} \dfrac{\partial g}{\partial x_j}.  \tag0-3$$

For semiclassical limits of quantum algebras, a natural version of the Gel'fand-Kirillov problem
is thus to ask whether the quotient field is isomorphic to a Poisson field of the form
$k_{\bfla}(x_1,\dots,x_n)$, or at least $K_{\bfla}(x_1,\dots,x_n)$ where $K$ is an extension field
of $k$. We establish the former for large classes of Poisson polynomial algebras, and the latter
for Poisson prime factors of these algebras.

In the last section of the paper, we introduce some invariants for Poisson fields,
with which we show that $k_{\bfla}(x_1,\dots,x_n)$ is never isomorphic to a
Poisson-Weyl field, and with which we can separate isomorphism classes of the
Poisson fields $k_{\bfla}(x_1,\dots,x_n)$ in many cases.
\enddefinition

\head 1. A finiteness theorem for torus-stable Poisson primes \endhead

In this section, we prove our finiteness theorem for the number of Poisson prime
ideals stable under a suitable torus action on an iterated
Poisson polynomial algebra. We begin by recalling the concept of a
Poisson polynomial algebra as introduced by Oh \cite{\Ohtwo}.

\definition{1.1\. Poisson polynomial algebras} Let $B$ be a Poisson
algebra. A {\it Poisson derivation\/} on $B$ is a ($k$-linear) map
$\alpha$ on $B$ which is a derivation with respect to both the
multiplication and the Poisson bracket, that is, $\alpha(ab)= \alpha(a)b+
a\alpha(b)$ and $\alpha(\{a,b\})= \{\alpha(a),b\}+ \{a,\alpha(b)\}$ for
$a,b\in B$. Suppose that $\delta$ is a derivation on $B$ such that
$$\delta(\{a,b\})= \{\delta(a),b\}+ \{a,\delta(b)\}+ \alpha(a)\delta(b)-
\delta(a)\alpha(b)  \tag 1-1$$
for $a,b\in B$. By \cite{\Ohtwo, Theorem 1.1} (after replacing our $B$
and $\alpha$ with $A$ and $-\alpha$), the Poisson structure on $B$ extends
uniquely to a Poisson algebra structure on the polynomial ring $A= B[x]$
such that
$$\{x,b\}= \alpha(b)x+ \delta(b)  \tag 1-2$$
 for $b\in B$. We write $A=
B[x;\alpha,\delta]_p$ to denote this situation, and we refer to $A$ as a
{\it Poisson polynomial algebra\/}.

The Poisson structure on $A$ extends uniquely to the Laurent polynomial
ring $B[x^{\pm1}]$, and is again determined by $\alpha$ and $\delta$. Hence, we
write $B[x^{\pm1};\alpha,\delta]_p$ for the ring $B[x^{\pm1}]$ equipped with this
structure, and we refer to it as a {\it Poisson Laurent polynomial
algebra\/}.

In either of the above cases, we omit $\delta$ from the notation if it is zero,
that is, we write $B[x;\alpha]_p$ and $B[x^{\pm1};\alpha]_p$ for
$B[x;\alpha,0]_p$ and $B[x^{\pm1};\alpha,0]_p$ respectively.

We will also need the converse part of \cite{\Ohtwo, Theorem 1.1}: if a
polynomial ring $A=B[x]$ supports a Poisson bracket such that $B$ is a
Poisson subalgebra and $\{x,B\} \subseteq Bx+B$, then $A=
B[x;\alpha,\delta]_p$ for suitable $\alpha$ and $\delta$.
\enddefinition

\proclaim{1.2\. Lemma} Let $A= B[x^{\pm1};\alpha]_p$ be a Poisson Laurent
polynomial algebra, and assume that $\alpha$ extends to a derivation $\alphahat$
on $A$ such that
$\alphahat(x)= sx$ for some nonzero $s\in k$. Then every $\alphahat$-stable
Poisson prime of $A$ is induced from a Poisson prime of $B$. \endproclaim

\demo{Proof} Let $P$ be an $\alphahat$-stable Poisson prime of $A$, and note that
$P\cap B$ is a Poisson prime of $B$. Then $(P\cap B)[x^{\pm1}]$ is an
$\alphahat$-stable Poisson prime of $A$, and we may pass to $A/(P\cap
B)[x^{\pm1}]$. Thus, without loss of generality, $P\cap B=0$, and we must show
that $P=0$.

If $P\ne 0$, then $P\cap B[x] \ne 0$. Choose a nonzero polynomial $p\in P\cap
B[x]$ of minimal degree, say $p= b_nx^n+ b_{n-1}x^{n-1}+ \cdots+ b_1x+b_0$ with
the $b_i\in B$ and $b_n\ne 0$. Note that $n>0$, because $P\cap B=0$. Now $P$
contains the polynomial
$$\alphahat(p)- \{x,p\}x^{-1}= nsb_nx^n+ (n-1)sb_{n-1}x^{n-1}+ \cdots+ sb_1x,$$
and hence also the polynomial
$$nsp- \bigl( \alphahat(p)- \{x,p\}x^{-1} \bigr) = sb_{n-1}x^{n-1}+ \cdots+
(n-1)sb_1x +nsb_0.$$
The latter must vanish, due to the minimality of $n$, and so $b_i=0$
for $i<n$. But then $b_nx^n= p\in P$ and so $b_n\in P$, contradicting the
assumption that
$P\cap B=0$. Therefore $P=0$, as required. \qed\enddemo

\proclaim{1.3\. Proposition} Let $A= B[x;\alpha,\delta]_p$ be a Poisson
polynomial algebra, and assume that $\alpha$ extends to a derivation $\alphahat$
on $A$ such that
$\alphahat(x)= sx$ for some nonzero $s\in k$. For each Poisson prime $Q$ of $B$,
there are at most two $\alphahat$-stable Poisson primes of $A$ that contract to
$Q$.
\endproclaim

\demo{Proof} Assume there exists an $\alphahat$-stable Poisson prime $P$ in $A$
that contracts to $Q$. For $b\in Q$, we have $\{x,b\}\in P$ and $\alpha(b)=
\alphahat(b)\in P$, whence $\delta(b)= \{x,b\}-\alpha(b)x\in P$, and so
$\alpha(b),\delta(b)\in Q$. It follows that $\{x,Q[x]\} \subseteq Q[x]$, from
which we see that $Q[x]$ is an $\alphahat$-stable Poisson prime of $A$. Hence, we
may pass to $A/Q[x]$ and then localize $B/Q$ to its quotient field. Thus, without
loss of generality,
$B$ is a field, and we must show that $A$ has at most two $\alphahat$-stable
Poisson primes.

Assume there exists a nonzero $\alphahat$-stable Poisson prime $P$ in $A$. Let
$n$ be the minimum degree of nonzero elements of $P$, and choose a monic
polynomial $p\in P$ of degree $n$, say  $p= x^n+ b_{n-1}x^{n-1}+ \cdots+
b_1x+b_0$ with the $b_i\in B$. Now $P$ contains the polynomial
$$\alphahat(p)-nsp = [\alpha(b_{n-1})-sb_{n-1}]x^{n-1} +\cdots+
[\alpha(b_1)-(n-1)sb_1]x + [\alpha(b_0)-nsb_0],$$
which must be zero by the minimality of $n$, and so $\alpha(b_{n-1})= sb_{n-1}$.
For any $b\in B$, the following polynomial lies in $P$:
$$\multline
\{p,b\}- n\alpha(b)p= \bigl[ n\delta(b)+ \{b_{n-1},b\} +(n-1)b_{n-1}\alpha(b)-
n\alpha(b)b_{n-1} \bigr] x^{n-1}  \\
 + \bigl[\text{lower terms}\bigr].  \endmultline$$
This polynomial must be zero, and so $n\delta(b)+ \{b_{n-1},b\}=
\alpha(b)b_{n-1}$. Thus, the element $d := \frac1n b_{n-1} \in B$ satisfies
$\alpha(d)=sd$ and $\{d,b\}= \alpha(b)d- \delta(b)$ for $b\in B$.

Set $y:= x+d$. Then $A$ is a Poisson polynomial algebra of the form $A=
B[y;\alpha]_p$. Further,
$\alphahat(y)= sy$, and so Lemma 1.2 implies that the only $\alphahat$-stable
Poisson prime of
$B[y^{\pm1};\alpha]_p$ is zero. Therefore, the only $\alphahat$-stable Poisson
primes of $A$ are
$\langle0\rangle$ and $\langle y\rangle$. \qed\enddemo

Our finiteness theorem parallels a
corresponding finiteness result of Letzter and the first author \cite{\specstrat,
Theorem 4.7}, which applies to torus actions on iterated skew polynomial algebras.
A key hypothesis in the latter theorem is that the automorphisms involved in the
skew polynomial structure must be restrictions of elements of the acting torus.
In the Poisson case, the corresponding ingredients are Poisson derivations, and
the relevant hypothesis relates these to the differential of the torus action. We
next recall the key facts about differentials of actions.

\definition{1.4\. The differential of a group action} Let $A$ be a $k$-algebra and
$G$ an algebraic group over
$k$, with Lie algebra $\gfrak$. Let $\alpha: G\rightarrow \Aut_{k\text{-alg}}(A)$
be a rational action of
$G$ on a $k$-algebra $A$ by $k$-algebra automorphisms. Thus,
$A$ is a directed union of finite dimensional $G$-stable subspaces $V_i$ such
that the induced maps
$\alpha_i :G\rightarrow GL(V_i)$ are morphisms of algebraic groups. In this
situation, the following hold:
\roster
\item The differentials $d\alpha_i:
\gfrak \rightarrow \glfrak(V_i)$ are compatible with inclusions $V_i\subseteq
V_j$, and they induce an action $d\alpha: \gfrak \rightarrow
\Der_k(A)$.
\item If $G$ is connected, the
$G$-stable subspaces of $A$ coincide with the
$\gfrak$-stable subspaces.
\endroster
That the $d\alpha_i$ are compatible with inclusions $V_i\subseteq V_j$ is a
routine check, as in \cite{\TauYu, Proposition 23.4.17}. One thus obtains a
linear action
$d\alpha: \gfrak \rightarrow \End_k(A)$, called the {\it differential of the
$G$-action\/}. Statement (2) is proved in \cite{\TauYu, Corollary 24.3.3} for the
case that $k$ is algebraically closed, but the latter hypothesis is not required.
The remainder of statement (1) is standard, but we have not located a precise
reference. It can be quickly obtained from two results in \cite{\TauYu}, as
follows.  For each $i$, the multiplication map
$V_i \otimes V_i \rightarrow V_i^2 \subseteq A$ is $G$-equivariant, and
so it is $\gfrak$-equivariant \cite{\TauYu, Proposition 23.4.17}.
Since the $G$-action on $V_i \otimes V_i$ is the diagonal one, so is the
$\gfrak$-action
\cite{\TauYu, Proposition 23.4.12}, from which we conclude that $\gfrak$ acts on
$A$ by derivations.

If $G$ is a torus, rationality of the action means that $A$ is the direct sum of its
$G$-eigenspaces, and the corresponding eigenvalues are rational characters of $G$ (e.g., see
\cite{\BrGo, Theorem II.2.7}). In this case, $A$ is also the direct sum of its
$\gfrak$-eigenspaces, and we have the following explicit description of the $\gfrak$-action. We
replace $G$ and $\gfrak$ by $H$ and $\hfrak$ to match our later notation, and we write $A_x$ for
the $x$-eigenspace of $A$, where $x$ is a member of the character group $X(H)$. Finally, we use
$(-|-)$ to denote the Euclidean inner product (or ``dot product'') in any $k^r$.
\enddefinition

\proclaim{1.5\. Lemma} Let $A$ be a $k$-algebra, equipped with a rational action
of a torus $H= (\kx)^r$ {\rm(}by $k$-algebra automorphisms\/{\rm)}. Identify
$\hfrak=\Lie H$ with $k^r$, and let $\hfrak$ act on $A$ by the differential of the
$H$-action. Further, identify $\ZZ^r$ with $X(H)$ via the natural pairing
$$\align \ZZ^r\times (\kx)^r &\longrightarrow \kx  \\
(m_1,\dots,m_r,h_1,\dots,h_r) &\longmapsto h_1^{m_1}h_2^{m_2}\cdots h_r^{m_r}.
\endalign$$
Then $\eta.a= (\eta|x)a$ for $\eta\in\hfrak$, $a\in A_x$, and $x\in X(H)$.

In particular, it follows that the $\hfrak$-action on $A$ commutes with the
$H$-action.
\qed \endproclaim

Readers who do not wish to delve into the full theory of actions of algebraic
groups may take the formula in Lemma 1.5 as the definition of the $\hfrak$-action
on $A$.

Whenever we have a torus $H$ acting rationally on a $k$-algebra $A$, we will
assume that its Lie algebra $\hfrak$ correspondingly acts on $A$ by the
differential of the $H$-action. We label the action of $H$ on $A$ a {\it rational
Poisson action\/} in case
$A$ is a Poisson algebra and $H$ acts rationally on $A$ by Poisson automorphisms.

\proclaim{1.6\. Lemma} Let $A$ be a Poisson algebra, equipped with a rational
Poisson action of a torus $H$. Then
$\hfrak=
\Lie H$ acts on $A$ by Poisson derivations. \endproclaim

\demo{Proof} Let $\eta\in \hfrak$, and let $a\in A_x$ and $b\in A_y$ for some
$x,y\in X(H)$. Since $h.\{a,b\}= \{h.a,h.b\}= x(h)y(h)\{a,b\}$ for all $h\in
H$, we have $\{a,b\}\in A_{x+y}$. Taking account of the  identifications in Lemma
1.5, we see that
$$\eta.\{a,b\}= (\eta|x+y)\{a,b\}= \bigl( (\eta|x)+ (\eta|y) \bigr) \{a,b\}=
\{\eta.a,b\}+ \{a,\eta.b\}.$$ Therefore, since $A$ is $X(H)$-graded, we conclude that $\eta.(-)$
is a Poisson derivation on $A$. \qed\enddemo

\proclaim{1.7\. Theorem} Let $A= k[x_1] [x_2;\alpha_2,\delta_2]_p \cdots
[x_n;\alpha_n,\delta_n]_p$ be an iterated Poisson polynomial algebra, supporting
a rational action by a torus $H$ such that $x_1,\dots,x_n$ are
$H$-eigen\-vec\-tors. Assume that there exist $\eta_1,\dots,\eta_n\in \hfrak= \Lie
H$ such that $\eta_i.x_j= \alpha_i(x_j)$ for $i>j$ and the $\eta_i$-eigenvalue of
$x_i$ is nonzero for each $i$. Then $A$ has at most $2^n$ $H$-stable
Poisson primes.
\endproclaim

\remark{Remark} Here the elements of $H$ are only assumed to act on $A$ by
$k$-algebra automorphisms, not necessarily by Poisson automorphisms. However, the
assumption of a Poisson action is needed in Corollary 1.8.
\endremark

\demo{Proof} Set $A_i= k[x_1] [x_2;\alpha_2,\delta_2]_p \cdots
[x_i;\alpha_i,\delta_i]_p$ for $i=0,1,\dots,n$. In view of \S1.4, the $H$-stable
Poisson primes in $A_i$ coincide with
the $\hfrak$-stable Poisson primes. Obviously $A_0=k$ has only one
$\hfrak$-stable Poisson prime. Now let $i<n$ and assume that $A_i$ has a finite
number, say $n_i$, of $\hfrak$-stable Poisson primes. It follows from the
relations $\eta_i.x_j=
\alpha_i(x_j)$ that the action of $\eta_i$ on $A_i$ coincides with $\alpha_i$.
Proposition 1.3 now implies that the number of $\hfrak$-stable
Poisson primes in $A_{i+1}$ is at most $2n_i$. The theorem follows. \qed\enddemo

\proclaim{1.8\. Corollary} Let $A= k[x_1] [x_2;\alpha_2,\delta_2]_p \cdots
[x_n;\alpha_n,\delta_n]_p$ be an iterated Poisson polynomial algebra, supporting
a rational Poisson action by a torus $H$ such that $x_1,\dots,x_n$ are
$H$-eigenvectors. Assume that there exist $\eta_1,\dots,\eta_n\in \hfrak= \Lie H$
such that $\eta_i.x_j= \alpha_i(x_j)$ for $i>j$ and the $\eta_i$-eigenvalue of
$x_i$ is nonzero for each $i$. Then $A$ satisfies the Poisson Dixmier-Moeglin
equivalence. \endproclaim

\demo{Proof} Theorem 1.7 and \cite{\pdixmo, Theorem 4.3}. \qed\enddemo

Theorem 1.7 and Corollary 1.8 can be extended to certain non-polynomial affine
Poisson algebras as follows.

\proclaim{1.9\. Proposition} Let $A$ be a Poisson algebra which is generated
{\rm(}as an algebra\/{\rm)} by a Poisson subalgebra $B$ together with a single
element $x$. Assume that $A$ supports a rational Poisson action by a torus $H$
such that
$B$ is $H$-stable and $x$ is an $H$-eigenvector. Moreover, assume that there is
some $\eta_0\in \hfrak= \Lie H$ such that $\{x,b\}- (\eta_0.b)x \in B$ for all
$b\in B$, and such that the $\eta_0$-eigenvalue of $x$ is nonzero. Then there are
at most twice as many $H$-stable Poisson primes in $A$ as in $B$.  \endproclaim

\demo{Proof} We show that $A$ is an epimorphic image of a Poisson polynomial ring
$\Ahat= B[X;\alpha,\delta]_p$ to which Proposition 1.3 applies. Let $\alpha$
denote the restriction of $\eta_0.(-)$ to $B$. Then, by Lemma 1.6, $\alpha$ is a
Poisson derivation on $B$, and, by hypothesis, $\delta(b) :=
\{x,b\}-
\alpha(b)x \in B$ for all $b\in B$. Since $\{x,-\}$ and $\alpha$ are
derivations, so is
$\delta$. For $b,b'\in B$, we compute that
$$\align \delta \bigl( \{b,b'\} \bigr) &= \{x,\{b,b'\}\}- \alpha(\{b,b'\})x  \\
 &= -\{b,\{b',x\}\}- \{b',\{x,b\}\}- \bigl( \{\alpha(b),b'\}+ \{b,\alpha(b')\}
\bigr) x  \\
 &= \bigl( \{\{x,b\},b'\}- \{\alpha(b)x,b'\} \bigr) + \{\alpha(b)x,b'\}-
\{\alpha(b),b'\}x  \\
 &\qquad\qquad  + \bigl( \{b,\{x,b'\}\}- \{b,\alpha(b')x\} \bigr) +
\{b,\alpha(b')x\}- \{b,\alpha(b')\}x  \\
 &= \{\delta(b),b'\}+ \{b,\delta(b')\}+ \alpha(b)\{x,b'\}+ \alpha(b')\{b,x\}  \\
 &= \{\delta(b),b'\}+ \{b,\delta(b')\}+ \alpha(b) \bigl( \alpha(b')x +\delta(b')
\bigr) +\alpha(b') \bigl( -\alpha(b)x -\delta(b) \bigr)  \\
 &= \{\delta(b),b'\}+ \{b,\delta(b')\}+ \alpha(b)\delta(b')- \delta(b)\alpha(b').
\endalign$$
Thus, the conditions for the existence of the Poisson polynomial ring
$\Ahat= B[X;\alpha,\delta]_p$ are verified.

Let $f\in X(H)$ be the $H$-eigenvalue of $x$. The action of $H$ on $B$ extends to
a rational action of $H$ on $\Ahat$ (by algebra automorphisms, at least) such
that $X$ is an $H$-eigenvector with $H$-eigenvalue $f$. (It is easily checked
that $H$ acts on $\Ahat$ by Poisson automorphisms, but we shall not need this
fact.) Since $x$ and $X$ have the same $H$-eigenvalue, they have the same
$\hfrak$-eigenvalue, and hence the same $\eta_0$-eigenvalue. Thus, the
$\eta_0$-eigenvalue of $X$ is nonzero. Since $\eta_0$ acts as a derivation on
$\Ahat$ extending $\alpha$, Proposition 1.3 now implies that for each Poisson
prime $Q$ of $B$, there are at most two $\eta_0$-stable Poisson primes of $\Ahat$
that contract to $Q$. All $H$-stable ideals of $\Ahat$ are $\hfrak$-stable and
thus $\eta_0$-stable, and so we conclude that for each $H$-stable Poisson prime
$Q$ of $B$, there are at most two $H$-stable Poisson primes of $\Ahat$ that
contract to $Q$. Thus, there are at most twice as many $H$-stable Poisson primes
in $\Ahat$ as in $B$.

Finally, the identity map on $B$ extends to a $k$-algebra surjection $\pi: \Ahat
\rightarrow A$ such that $\pi(X)= x$. Obviously $\pi$ preserves brackets of
elements of $B$, and
$$\pi(\{X,b\})= \alpha(b)x+\delta(b)= \{x,b\}= \{\pi(X),\pi(b)\}$$
for $b\in B$, from which we see that $\pi$ is a Poisson homomorphism. By
construction, $\pi$ is also $H$-equivariant. Hence, the set map $\pi^{-1}$ embeds
the collection of $H$-stable Poisson primes of $A$ into the collection of
$H$-stable Poisson primes of $\Ahat$. Therefore, there are
at most twice as many $H$-stable Poisson primes in $A$ as in $B$. \qed\enddemo

\proclaim{1.10\. Theorem} Let $A$ be a Poisson algebra, equipped with a rational
Poisson action by a torus $H$. Assume that $A$ is generated by $H$-eigenvectors
$x_1,\dots,x_n$, and that there exist $\eta_1,\dots,\eta_n \in \hfrak= \Lie H$
such that
\roster
\item $\{x_i,x_j\}- (\eta_i.x_j)x_i \in k\langle x_1,\dots,x_{i-1} \rangle$
for all $i>j$;
\item For all $i$, the $\eta_i$-eigenvalue of $x_i$ is nonzero.
\endroster
Then $A$ has at most $2^n$ $H$-stable Poisson primes, and $A$ satisfies the
Poisson Dixmier-Moeglin equivalence. \endproclaim

\demo{Proof} The first conclusion is clear when $n=0$. Now let $n>0$, and assume
that the subalgebra $B:= k\langle x_1,\dots,x_{n-1}\rangle$ has at most $2^{n-1}$
$H$-stable Poisson primes. Note that the map $\delta_n := \{x_n,-\}-
(\eta_n.-)x_n$ is a derivation on $A$. Since, by hypothesis,
$\delta_n(x_j) \in B$ for all $j<n$, it follows that $\delta_n(B) \subseteq B$.
Since the
$\eta_n$-eigenvalue of $x_n$ is nonzero, Proposition 1.9 implies that $A$ has at
most twice as many $H$-stable Poisson primes as $B$, thus at most $2^n$.

The final conclusion now follows from \cite{\pdixmo, Theorem 4.3}. \qed
\enddemo

\head 2. Poisson polynomial algebras satisfying the Poisson Dixmier-Moeglin
equivalence \endhead

We show that semiclassical limits of many standard quantum algebra constructions
yield Poisson polynomial algebras to which Theorem 1.7 and Corollary 1.8 apply.

\definition{2.1\. Semiclassical limits} Suppose that $R$ is a commutative
principal ideal domain, containing $k$, and that $\hsl\in R$ with $\hsl R$ a maximal ideal of
$R$. If $B$ is a torsionfree $R$-algebra for which the
quotient $A:= B/\hsl B$ is commutative, then there is a well-defined bilinear map
$\frac{1}{\hsl} [-,-]: B\times B \rightarrow B$, which induces a Poisson bracket
on $A$ (e.g., see \cite{\BrGo, \S III.5.4}). The Poisson algebra $A$ is known as
the {\it semiclassical limit\/} (or {\it quasiclassical limit\/}) of $B$, or of
the family of algebras $(B/\mfrak B)_{\mfrak\in\max R}$.

There are two standard choices for $R$ in quantum algebra constructions. In
single parameter cases, we take $R$ to be a Laurent polynomial ring $k[q,q^{-1}]$
and $\hsl= q-1$, while multiparameter cases are usually best handled by taking $R$
to be a formal power series algebra
$\khh$. In the latter situation, we use the abbreviation
$$e(\alpha) := \operatorname{exp}(\alpha\hsl)= \sum_{i=0}^\infty \frac{1}{i!}
\alpha^i \hsl^i$$
for $\alpha\in k$. Note that $e(\alpha+\beta)= e(\alpha)
e(\beta)$ for $\alpha,\beta\in k$.
\enddefinition

\definition{2.2\. Semiclassical limits of quantum affine spaces} {\bf (a)}
Suppose that $\bfq= (q_{ij})$ is an $n\times n$
multiplicatively antisymmetric matrix over $k$, that is,
$q_{ii}=1$ and $q_{ji}=q_{ij}^{-1}$ for all $i$, $j$. The corresponding
multiparameter quantized coordinate ring of affine $n$-space is the $k$-algebra
$\O_{\bfq}(k^n)$ with generators $x_1,\dots,x_n$ and relations
$x_ix_j= q_{ij}x_jx_i$ for all $i$, $j$. Similarly, if $\bfq$ is a
multiplicatively antisymmetric matrix over a commutative ring $R$, we can form
the $R$-algebra $\O_{\bfq}(R^n)$. Observe that $\O_{\bfq}(R^n)$ is an iterated
skew polynomial algebra over $R$, and hence a free $R$-module.

{\bf (b)} To write the semiclassical limits of the above
algebras, a change of notation is convenient. Let $\bfq$ now be an
(additively) antisymmetric matrix in $M_n(k)$. Since the matrix $e(\bfq)=
(e(q_{ij}))$ is a multiplicatively antisymmetric matrix over $\khh$, we can form
the $\khh$-algebra $B=
\O_{e(\bfq)}(\khh^n)$. As noted in (a), $B$ is a free $\khh$-module, and hence it
is torsionfree over $\khh$. We identify the quotient
$A= B/\hsl B$ with the polynomial algebra $k[x_1,\dots,x_n]$. Since this
algebra is commutative, it inherits a Poisson bracket such that $\{x_i,x_j\}=
q_{ij}x_ix_j$ for all $i$,
$j$.

{\bf (c)} There is a rational action of the torus $H= (\kx)^n$ on $A$ such that
$$(h_1,\dots,h_n).x_i= h_ix_i$$
for $(h_1,\dots,h_n) \in H$ and $i=1,\dots,n$. This action preserves the Poisson
bracket on the indeterminates, that is, $h.\{x_i,x_j\}= \{h.x_i,h.x_j\}$ for
$h\in H$ and all $i$, $j$. Consequently, it is a Poisson
action. In this case, Theorem 1.7 is not needed, since $A$ clearly has
exactly $2^n$ $H$-stable primes, namely the ideals $\langle x_i \mid i\in
I\rangle$ for $I \subseteq \{1,\dots,n\}$. That
$A$ satisfies the Poisson Dixmier-Moeglin equivalence was shown in \cite{\pdixmo,
Example 4.6}.
\enddefinition

\definition{2.3\. Semiclassical limits of quantum matrices} {\bf
(a)} Given a nonzero scalar $\lambda\in \kx$ and a multiplicatively antisymmetric
matrix $\bfp= (p_{ij})
\in M_n(\kx)$, the multiparameter quantum $n\times n$ matrix
algebra $\OlpMn$ is the $k$-algebra with generators $X_{ij}$ for
$i,j=1,\dots,n$ and relations
$$X_{lm}X_{ij} = \cases p_{li}p_{jm}X_{ij}X_{lm} +
(\lambda -1)p_{li}X_{im}X_{lj} &\quad (l>i,\ m>j)\\
\lambda p_{li}p_{jm}X_{ij}X_{lm} &\quad (l>i,\ m\le j)\\
p_{jm}X_{ij}X_{lm} &\quad (l=i,\ m>j).\endcases  \tag 2-1$$
The standard single parameter case is recovered when $\lambda= q^{-2}$ and
$p_{ij}=q$ for all $i>j$. When $\lambda=1$, we just have a multiparameter
quantum affine $n^2$-space, $\O_{\bfq}(k^{n^2})$, for suitable $\bfq$.

{\bf (b)} Now let $\bfp$ be an antisymmetric matrix in $M_n(k)$,
and $\lambda\in k$ an arbitrary scalar. Form the algebra $B=
\O_{e(\lambda),e(\bfp)}(M_n(\khh))$, and identify the quotient $A= B/\hsl B$ with
the polynomial algebra over $k$ in the indeterminates $X_{ij}$,
that is, $A=\OMn$. One can check directly that $B$ is an iterated skew polynomial
algebra over $\khh$. Alternatively, it is known that
$\O_{e(\lambda),e(\bfp)}(M_n(k((\hsl))))$ is an iterated skew polynomial algebra
over the field $k((\hsl))$, and one observes that the automorphisms and skew
derivations of this structure map the relevant $\khh$-subalgebras into
themselves. Either way, we conclude that $B$ is torsionfree over $\khh$.

Now $\OMn$ inherits a Poisson bracket
such that
$$\{ X_{lm}, X_{ij} \}= \cases (p_{li}+p_{jm})X_{ij}X_{lm} +
\lambda X_{im}X_{lj} &\quad (l>i,\ m>j)\\
(\lambda+ p_{li}+p_{jm})X_{ij}X_{lm} &\quad (l>i,\ m\le j)\\
p_{jm}X_{ij}X_{lm} &\quad (l=i,\ m>j).\endcases  \tag 2-2$$
When $\lambda=0$, we have a semiclassical limit of a quantum affine $n^2$-space, a
case covered in \S2.2. Hence, we now assume that $\lambda\ne 0$.

Observe that
$$\multline
\{ X_{lm}\,,\, k[X_{ij} \mid (i,j) \lexle (l,m)] \} \subseteq \\
 k[X_{ij} \mid (i,j)
 \lexle (l,m)] X_{lm}+ k[X_{ij} \mid (i,j) \lexle (l,m)] \endmultline$$  
 for all $l$, $m$, and so when the $X_{ij}$ are adjoined in
lexicographic order, $\OMn$ is an iterated Poisson polynomial algebra of the form
$$\OMn= k[X_{11}][X_{12};\alpha_{12},\delta_{12}]_p \cdots [X_{nn}; \alpha_{nn},
\delta_{nn}]_p \,.$$
In view of (2-2), we have
$$\alpha_{lm}(X_{ij})= \cases (p_{li}+p_{jm})X_{ij} &\quad (l>i,\ m>j)\\
(\lambda+ p_{li}+p_{jm})X_{ij} &\quad (l>i,\ m\le j)\\
p_{jm}X_{ij} &\quad (l=i,\ m>j).\endcases  \tag2-3$$

{\bf (c)} There is a rational action of the torus $H= (\kx)^{2n}$ on $\OMn$ such
that
$$h.X_{ij}= h_i h_{n+j}
X_{ij} \qquad\qquad (\, h=(h_1,\dots,h_{2n}) \in H\,)$$ for all $i$, $j$, and it is clear from
(2-2) that this is a Poisson action. If we identify $\ZZ^{2n}$ with $X(H)$ as in Lemma 1.5, then
each $X_{ij}$ has $H$-eigenvalue $\eps_i+\eps_{n+j}$, where $(\eps_1,\dots,\eps_{2n})$ is the
canonical basis for $\ZZ^{2n}$. Hence, the differential of the $H$-action gives an action of
$\hfrak= \Lie H= k^{2n}$ on $\OMn$ by derivations such that
$$\eta.X_{ij}= (\eta_i+ \eta_{n+j})
X_{ij} \qquad\qquad (\, \eta= (\eta_1,\dots,\eta_{2n}) \in \hfrak\,)$$
for all $i$, $j$.

{\bf (d)} For $l,m= 1,\dots,n$, define $\eta_{lm} \in \hfrak$ as follows:
$$\eta_{lm}= (p_{l1},\dots,p_{ln}, p_{1m},\dots,p_{m-1,m}, \lambda,
\lambda+p_{m+1,m}, \dots, \lambda+p_{nm}).$$ Then observe that $\eta_{lm}.X_{ij}=
\alpha_{lm}(X_{ij})$ for $(i,j) \lexle (l,m)$, and that the $\eta_{lm}$-eigenvalue of $X_{lm}$ is
$\lambda$. Since we have assumed that $\lambda\ne 0$, the hypotheses of Theorem 1.7 and Corollary
1.8 are satisfied. Therefore $A$ has at most $2^{n^2}$ $H$-stable Poisson primes, and $A$
satisfies the Poisson Dixmier-Moeglin equivalence. The case of the latter result with $n=2$ and
the standard Poisson bracket was established by Oh in \cite{\Oh, Proposition 2.3}.
\enddefinition

\definition{2.4\. Semiclassical limits of quantum symplectic and
even-di\-men\-sion\-al euclidean spaces} Multiparameter versions of the mentioned
quantum algebras are instances of the algebras $K^{P,Q}_{n, \Gamma}(k)$ introduced
by Horton
\cite{\Hrt}, and we treat that general class.

{\bf (a)} Let $\Gamma = (\gamma_{ij}) \in M_n(\kx)$ be a
multiplicatively antisymmetric matrix, and let $P=(p_1, \dots p_n)$ and $Q =
(q_1, \dots q_n)$ be vectors in $(\kx)^n$ such that $p_i \neq q_i$ for all $i$.
Then $K^{P,Q}_{n, \Gamma}(k)$ is the
$k$-algebra with generators $x_1, y_1, \dots, x_n, y_n$ and relations
$$\alignedat2 y_iy_j & = \gamma_{ij}y_jy_i &\qquad\qquad\qquad&(\text{all\ } i,j)
\\ x_iy_j &= p_j\gamma_{ji}y_jx_i &&(i<j) \\
x_iy_j &= q_j \gamma_{ji}y_jx_i && (i>j) \\
x_ix_j &= q_ip_j^{-1}\gamma_{ij}x_jx_i && (i<j) \\
x_iy_i &= q_iy_ix_i + \sum _{\ell < i}(q_{\ell}-p_{\ell}) y_{\ell}x_{\ell} &&(\text{all\ }
i).\endalignedat  \tag 2-4$$
See \cite{\Hrt, Examples 1.3--1.7} for the choices of parameters
which yield the standard quantum symplectic and even-dimensional euclidean spaces, and related
algebras. This construction can be performed over a commutative ring $R$, assuming the $p_i$,
$q_i$, and $\gamma_{ij}$ are units in $R$, and as in \cite{\Hrt, Proposition 2.5}, $K^{P,Q}_{n,
\Gamma}(R)$ is an iterated skew polynomial algebra over $R$.

{\bf (b)} Now let $\Gamma$ be an antisymmetric matrix in $M_n(k)$, and let $P$
and $Q$ be vectors in $k^n$, with $p_i \neq q_i$ for all $i$. Form the algebra
$B= K^{e(P),e(Q)}_{n,e(\Gamma)}(\khh)$, and identify the quotient $A= B/\hsl B$
with the polynomial algebra $k[x_1,y_1,\dots,x_n,y_n]$. Now $B$ is an iterated
skew polynomial algebra over $\khh$, so it is a torsionfree $\khh$-module. Hence,
$A$ inherits a Poisson bracket such that
$$\alignedat2 \{y_i,y_j\} & = \gamma_{ij}y_iy_j &\qquad\qquad\qquad&(\text{all\
} i,j) \\
\{x_i,y_j\} &= (p_j+\gamma_{ji})x_iy_j &&(i<j) \\
\{x_i,y_j\} &= (q_j+ \gamma_{ji})x_iy_j && (i>j) \\
\{x_i,x_j\} &= (q_i-p_j+\gamma_{ij})x_ix_j && (i<j) \\
\{x_i,y_i\} &= q_ix_iy_i + \sum _{\ell < i}(q_{\ell}-p_{\ell}) x_{\ell}y_{\ell} &&(\text{all\ }
i).\endalignedat  \tag 2-5$$
This Poisson algebra $A$ was introduced by Oh in \cite{\Ohthree} and
denoted $A^{P,Q}_{n, \Gamma}(k)$. It is an iterated Poisson polynomial algebra of the form
$$A= k[x_1][y_1;\alpha_1,\delta_1]_p [x_2;\alpha'_2]_p[y_2;\alpha_2,\delta_2]_p
\cdots [x_n;\alpha'_n]_p[y_n;\alpha_n,\delta_n]_p\,,$$
such that
$$\alignedat2 \alpha_j(x_i) &= (-p_j+ \gamma_{ij})x_i  &\alpha_j(y_i) &=
\gamma_{ji}y_i  \\
\alpha_j(x_j) &= -q_jx_j  \\
\alpha'_j(x_i) &= (-q_i+p_j+\gamma_{ji})x_i &\qquad\qquad\qquad   \alpha'_j(y_i) &=
(q_i+\gamma_{ij})y_i
\endalignedat  \tag2-6$$
for all $i<j$.

{\bf (c)} There is a rational Poisson action of the torus $H= (\kx)^{n+1}$ on $A$
such that
$$\alignedat2 h.x_i &= h_ix_i  &\qquad\qquad\qquad h.y_i &= h_1h_{n+1}h_i^{-1}y_i
\endalignedat$$
for $h\in H$. Then $\hfrak= \Lie H= k^{n+1}$ acts on
$A$ by derivations such that
$$\alignedat2 \eta.x_i &= \eta_ix_i  &\qquad\qquad\qquad \eta.y_i &=
(\eta_1+\eta_{n+1}-\eta_i) y_i \endalignedat$$
for $\eta\in\hfrak$.

{\bf (d)} Define $\eta_j,\eta'_j \in\hfrak$ as follows:
$$\xalignat2  \eta_1 &= (-q_1,0,\dots,0,1)  \\
\eta_j &= (-p_j+ \gamma_{1j},\dots,-p_j+ \gamma_{j-1,j}, -q_j,
0,\dots,0, \gamma_{j1})  &&(j>1)  \\
\eta'_j &= (-q_1+p_j+\gamma_{j1},\dots,-q_n+p_j+\gamma_{jn},q_1+\gamma_{1j})
&&(j>1).
\endxalignat$$
Note first that $\eta_1.x_1= \alpha_1(x_1)$, and that the $\eta_1$-eigenvalue of $y_1$ is $1$. For
$j>1$, we have $\eta_j.x_i= \alpha_j(x_i)$ for $i\le j$ and $\eta_j.y_i= \alpha_j(y_i)$ for $i<j$,
and the $\eta_j$-eigenvalue of $y_j$ is $q_j-p_j$. Finally, we have $\eta'_j.x_i= \alpha'_j(x_i)$
and $\eta'_j.y_i= \alpha'_j(y_i)$ for $i<j$, and the $\eta'_j$-eigenvalue of $x_j$ is $p_j-q_j$.
Thus, the hypotheses of Theorem 1.7 and Corollary 1.8 are satisfied. We conclude that $A$ has at
most $2^{2n}$ $H$-stable Poisson primes, and that it satisfies the Poisson Dixmier-Moeglin
equivalence. The case $n=2$ of the latter result was established by Oh in \cite{\Ohtwo, Theorem
3.5}.
\enddefinition

\definition{2.5\. Semiclassical limits of quantum odd-dimensional euclidean
spaces} Mul\-ti\-par\-a\-meter versions of quantum euclidean spaces in the
odd-dimensional case can be constructed analogously to the even-dimensional case
treated in \cite{\Hrt}. Since these algebras have not (to our knowledge) appeared
in the literature, we take the opportunity to introduce them here.

{\bf (a)} As for the $2n$-dimensional case, let $\Gamma = (\gamma_{ij}) \in
M_n(\kx)$ be a multiplicatively antisymmetric matrix, and let $P=(p_1, \dots
p_n)$ and $Q = (q_1, \dots q_n)$ be vectors in $(\kx)^n$ such that $p_i \neq q_i$
for all $i$. Further, let $\lambda\in k$, and assume that each $p_i$ has a
square root in $\kx$, which we fix and label $p_i^{1/2}$. Define
$K^{P,Q,\lambda}_{n,\Gamma}(k)$ to be the $k$-algebra with generators
$z_0,x_1,y_1,\dots,x_n,y_n$ and relations
$$\alignedat2 z_0x_i &= p_i^{-1/2}x_iz_0 &\qquad\qquad\qquad&(\text{all\ } i)  \\
z_0y_i &= p_i^{1/2}y_iz_0 &\qquad\qquad\qquad&(\text{all\ } i)  \\
y_iy_j & = \gamma_{ij}y_jy_i &\qquad\qquad\qquad&(\text{all\ } i,j) \\
x_iy_j &= p_j\gamma_{ji}y_jx_i &&(i<j) \\
x_iy_j &= q_j \gamma_{ji}y_jx_i && (i>j) \\
x_ix_j &= q_ip_j^{-1}\gamma_{ij}x_jx_i && (i<j) \\
x_iy_i &= q_iy_ix_i + \sum _{\ell < i}(q_{\ell}-p_{\ell}) y_{\ell}x_{\ell}
+\lambda z_0^2  &&(\text{all\ } i).
\endalignedat  \tag 2-7$$
The standard single parameter algebra corresponds to the case where the $p_i= q^{-2}$, the
$q_i=1$, the $\gamma_{ij}= q^{-1}$ for $i<j$, and $\lambda= (q-1)q^{n-(1/2)}$. (This requires a
change of variables, as in \cite{\Mus, \S\S2.1, 2.2}, \cite{\Ohcat, Example 5}, or \cite{\Hrt,
Example 1.5}.) On the other hand, if we take $\lambda=1$ and all the $p_i=1$ (with $p_i^{1/2}=1$),
then $z_0$ is central in $K^{P,Q,\lambda}_{n,\Gamma}(k)$, and $K^{P,Q,\lambda}_{n,\Gamma}(k) /
\langle z_0-1\rangle$ is the multiparameter quantized Weyl algebra $A_n^{Q,\Gamma}(k)$ (see, e.g.,
\cite{\BrGo, \S I.2.6}).

This construction
can be performed over a commutative ring $R$, assuming the relevant parameters
are units in $R$ and the $p_i$ have square roots in $R$. As in the
even-dimensional case \cite{\Hrt, Proposition 2.5},
$K^{P,Q,\lambda}_{n,\Gamma}(R)$ is an iterated skew polynomial algebra over $R$.

{\bf (b)} Now let $\Gamma$ be an antisymmetric matrix in $M_n(k)$, let $P$
and $Q$ be vectors in $k^n$ with $p_i \neq q_i$ for all $i$, and let $\lambda\in
k$. Form the algebra
$B= K^{e(P),e(Q),\lambda\hsl}_{n,e(\Gamma)}(\khh)$, and identify the quotient $A=
B/\hsl B$ with the polynomial algebra $k[z_0,x_1,y_1,\dots,x_n,y_n]$. Here we
have used $\lambda\hsl$ rather than $e(\lambda)$ to ensure commutativity of
$B/\hsl B$, and we take $e(p_i/2)$ as the chosen square root of $e(p_i)$. Since
$B$ is an iterated skew polynomial algebra over $\khh$, it is a torsionfree
$\khh$-module. Hence, $A$ inherits a Poisson bracket such that
$$\alignedat2 \{z_0,x_i\} &= -(p_i/2)z_0x_i  &&(\text{all\ } i)  \\
\{z_0,y_i\} &= (p_i/2)z_0y_i  &&(\text{all\ } i)  \\
\{y_i,y_j\} & = \gamma_{ij}y_iy_j &\qquad\qquad\qquad&(\text{all\
} i,j) \\
\{x_i,y_j\} &= (p_j+\gamma_{ji})x_iy_j &&(i<j) \\
\{x_i,y_j\} &= (q_j+ \gamma_{ji})x_iy_j && (i>j) \\
\{x_i,x_j\} &= (q_i-p_j+\gamma_{ij})x_ix_j && (i<j) \\
\{x_i,y_i\} &= q_ix_iy_i + \sum _{\ell < i}(q_{\ell}-p_{\ell}) x_{\ell}y_{\ell} +\lambda z_0^2
&&(\text{all\ } i).\endalignedat  \tag 2-8$$
This Poisson algebra is an iterated Poisson
polynomial algebra of the form
$$A= k[z_0][x_1;\alpha'_1]_p[y_1;\alpha_1,\delta_1]_p
[x_2;\alpha'_2]_p[y_2;\alpha_2,\delta_2]_p
\cdots [x_n;\alpha'_n]_p[y_n;\alpha_n,\delta_n]_p\,.$$

{\bf (c)} There is a rational Poisson action of the torus $H= (\kx)^{n+1}$ on $A$
such that
$$\alignedat2 h.x_i &= h_ix_i  &\qquad\qquad\qquad h.z_0 &= h_{n+1}z_0  \\
h.y_i &= h_{n+1}^2h_i^{-1}y_i   \endalignedat$$
for $h\in H$. We leave it to the reader to check that the hypotheses of Theorem
1.7 and Corollary 1.8 are satisfied. We conclude that $A$ has at most $2^{2n+1}$
$H$-stable Poisson primes, and that it satisfies the Poisson Dixmier-Moeglin
equivalence.
\enddefinition

\definition{2.6\. Semiclassical limits of quantum symmetric matrices} Fix a
positive integer $n$.

{\bf (a)} Coordinate rings of quantum symmetric $n\times n$ matrices
have been introduced by Noumi \cite{\Nou, Theorem 4.3,
Proposition 4.4, and comments following the proof} and Kamita \cite{\Kam,
Theorem 0.2}. As in \cite{\GoYa, \S5.5}, we take the case of Noumi's algebra with
all parameters equal to $1$, which agrees with Kamita's algebra after
interchanging the scalar parameter $q$ with $q^{-1}$. This is a $k$-algebra with
generators  $y_{ij}$ for $1\le i\le j\le n$. If the construction is instead
performed over a rational function field $k(q)$, the $k[q^{\pm1}]$-subalgebra
$B$ generated by the $y_{ij}$ is then an iterated skew polynomial algebra over
$k[q^{\pm1}]$, and the quotient
$A= B/(q-1)B$ can be identified with the polynomial ring $k[y_{ij} \mid 1\le
i\le j\le n]$. Hence, $A$ inherits a Poisson bracket, which has the following
form, as calculated in \cite{\GoYa, \S5.5}:
$$\{y_{ij}, y_{lm}\} =
\bigl( \sign(l-j) + \sign(m-i) \bigr) y_{il} y_{jm}+ \bigl( \sign(l-i) + \sign(m-j) \bigr) y_{im}
y_{jl}  \tag2-9$$
for $i\le j$ and $l\le m$, where $\sign(t)$ is $1$, $0$, or $-1$ according as
$t$ is positive, zero, or negative, and where $y_{ts}= y_{st}$ if needed. The Poisson algebra $A$
is an iterated Poisson polynomial algebra of the form
$$A= k[y_{11}][y_{12};\alpha_{12},\delta_{12}]_p \cdots
[y_{nn};\alpha_{nn},\delta_{nn}]_p\,,$$ where the indeterminates $y_{ij}$ (for $1\le i\le j\le n$)
have been adjoined in lexicographic order, and where
$$\alpha_{lm}(y_{ij})= \cases  -y_{ij}  &\quad \bigl( (i=l<j<m)
\text{\ or\ } (i<l<j=m) \text{\ or\ } (i<j=l<m) \bigr)  \\
-2y_{ij}  &\quad \bigl( (i=j=l<m) \text{\ or\ } (i<j=l=m) \bigr)  \\
0  &\quad (\text{otherwise})  \endcases  \tag2-10$$
for $l\le m$ and $i\le j$ with $(i,j) \lexle
(l,m)$.

{\bf (b)} There is a rational Poisson action of the torus $H= (\kx)^n$ on $A$ such that $h.y_{ij}=
h_ih_jy_{ij}$ for all $h\in H$ and all $i$, $j$. Then $\hfrak= \Lie H= k^n$ acts on $A$ by
derivations such that $\eta.y_{ij}= (\eta_i+\eta_j)y_{ij}$ for all $\eta\in\hfrak$ and all $i$,
$j$. Let $(\eps_1,\dots,\eps_n)$ denote the canonical basis for $\hfrak$, and set $\eta_{lm}=
-\eps_l-\eps_m$ for $1\le l\le m\le n$. Then $\eta_{lm}.y_{ij}= \alpha_{lm}(y_{ij})$ for $(i,j)
\lexle (l,m)$, and the $\eta_{lm}$-eigenvalue of $y_{lm}$ is either $-2$ or $-4$ (depending on
whether $l<m$ or $l=m$). Thus, the hypotheses of Theorem 1.7 and Corollary 1.8 are satisfied. We
conclude that $A$ has at most $2^{n(n+1)/2}$ $H$-stable Poisson primes, and that it satisfies the
Poisson Dixmier-Moeglin equivalence.
\enddefinition

\definition{2.7\. Semiclassical limits of quantum antisymmetric matrices} Fix a
positive integer $n$.

 {\bf (a)} The coordinate ring of quantum antisymmetric $n\times n$ matrices
was introduced by Strickland in \cite{\St, Section 1}; it is a
$k$-algebra with generators $y_{ij}$ for $1\le i<j\le n$ and relations involving
a scalar $q\in \kx$. If the construction is instead performed over $k(q)$, the
$k[q^{\pm1}]$-subalgebra $B$ generated by the $y_{ij}$ is then an iterated skew
polynomial algebra over $k[q^{\pm1}]$, and the quotient
$A= B/(q-1)B$ can be identified with the polynomial ring $k[y_{ij} \mid 1\le
i<j\le n]$. Hence, $A$ inherits a Poisson bracket, which has the following form,
as noted in \cite{\GoYa, \S5.6(b)}. (The factor $2$ in \cite{op.~cit.} does not
appear here, because we are using Strickland's construction without changing $q$
to $q^{1/2}$.)
$$\{y_{ij}, y_{lm}\} =
\bigl( \sign(l-j) + \sign(m-i) \bigr) y_{il} y_{jm}- \bigl( \sign(l-i) + \sign(m-j) \bigr) y_{im}
y_{jl}  \tag2-11$$
for $i<j$ and $l<m$, where $y_{ts}= -y_{st}$ and $y_{ss}=0$ if needed. The
Poisson algebra $A$ is an iterated Poisson polynomial algebra of the form
$$A= k[y_{12}][y_{13};\alpha_{13},\delta_{13}]_p \cdots
[y_{n-1,n};\alpha_{n-1,n},\delta_{n-1,n}]_p\,,$$ where the indeterminates $y_{ij}$ (for $1\le i<
j\le n$) have been adjoined in lexicographic order, and where
$$\alpha_{lm}(y_{ij})= \cases  -y_{ij}  &\qquad (\text{if\ } {\bigm|} \{i,j\}\cap
\{l,m\} {\bigm|}= 1)  \\   0  &\qquad (\text{otherwise})  \endcases  \tag2-12$$
for $l<m$ and
$i<j$ with $(i,j) \lexle (l,m)$.

{\bf (b)} There is a rational Poisson action of the torus $H= (\kx)^n$ on $A$ such that $h.y_{ij}=
h_ih_jy_{ij}$ for all $h\in H$ and all $i$, $j$. Then $\hfrak= \Lie H= k^n$ acts on $A$ by
derivations such that $\eta.y_{ij}= (\eta_i+\eta_j)y_{ij}$ for all $\eta\in\hfrak$ and all $i$,
$j$. Let $(\eps_1,\dots,\eps_n)$ denote the canonical basis for $\hfrak$, and set $\eta_{lm}=
-\eps_l-\eps_m$ for $1\le l<m\le n$. Then $\eta_{lm}.y_{ij}= \alpha_{lm}(y_{ij})$ for $(i,j)
\lexle (l,m)$, and the $\eta_{lm}$-eigenvalue of $y_{lm}$ is $-2$. Thus, the hypotheses of Theorem
1.7 and Corollary 1.8 are satisfied. We conclude that $A$ has at most $2^{n(n-1)/2}$ $H$-stable
Poisson primes, and that it satisfies the Poisson Dixmier-Moeglin equivalence.
\enddefinition

\head 3. Fraction fields of Poisson prime quotients \endhead

We now turn to the Poisson structure of fields of fractions of Poisson prime
quotients of iterated Poisson polynomial rings.

\definition{3.1} Recall from \S0.4 the notation $k_{\bfla}[x_1,\dots,x_n]$ for the
Poisson algebra based on the polynomial ring $k[x_1,\dots,x_n]$ with $\{x_i,x_j\}=
\lambda_{ij}x_ix_j$ for all $i$, $j$, where $\bfla= (\lambda_{ij})$ is an antisymmetric matrix in
$M_n(k)$. The corresponding Poisson Laurent polynomial algebra and Poisson field are denoted
$k_\bfla[x_1^{\pm1},\dots,x_n^{\pm1}]$ and $k_\bfla(x_1,\dots,x_n)$, respectively.

The algebra $k_\bfla[x_1^{\pm1},\dots,x_n^{\pm1}]$ can be identified with the group algebra
$k\Gamma$, where $\Gamma= \ZZ^n$, by writing  monomials in the $x_i$ in the form $x^\alpha=
x_1^{\alpha_1}x_2^{\alpha_2} \cdots x_n^{\alpha_n}$ for elements $\alpha=
(\alpha_1,\dots,\alpha_n) \in \Gamma$. There is then an antisymmetric bilinear form $b:
\Gamma\times\Gamma \rightarrow k$ such that
$$b(\alpha,\beta)= \sum_{i,j=1}^n \alpha_i \lambda_{ij} \beta_j$$
for $\alpha,\beta\in \Gamma$, and $\{x^\alpha,x^\beta\}= b(\alpha,\beta) x^{\alpha+\beta}$ for
$\alpha,\beta \in\Gamma$. Conversely, if $b$ is any $k$-valued antisymmetric bilinear form on
$\Gamma$, there is a Poisson bracket on $k\Gamma$ such that $\{x^\alpha,x^\beta\}=
b(\alpha,\beta)x^{\alpha+\beta}$ for $\alpha,\beta \in \Gamma$. We denote this Poisson algebra by
$k_b\Gamma$. The following facts about $k_b\Gamma$ are well known. See, for instance, \cite{\Van,
Lemma 1.2} where they are proved in the case $k=\CC$; the arguments are valid over arbitrary base
fields.
\enddefinition

\proclaim{3.2\. Lemma} Let  $\Gamma= \ZZ^n$ for some $n\in\NN$, let $b$ be a
$k$-valued antisymmetric bilinear form on $\Gamma$, and let $k_b\Gamma$ be the
Poisson algebra based on $k\Gamma$ described in {\rm \S3.1}. Set
$$\Gamma_b= \{\alpha\in \Gamma\mid b(\alpha,-) \equiv 0\},$$
a subgroup of $\Gamma$. Then $Z_p(k_b\Gamma)= k\Gamma_b$, and every Poisson ideal
of $k_b\Gamma$ is generated by its intersection with $Z_p(k_b\Gamma)$.
\qed\endproclaim

\proclaim{3.3\. Corollary} Let $\Gamma$ and $b$ be as in Lemma {\rm 3.2}. If
$Z_p(k_b\Gamma)= k$, then $k_b\Gamma$ is {\rm Poisson
simple}, that is, its only Poisson ideals are $0$ and $k_b\Gamma$.
\qed\endproclaim

We can now determine the structure of the fields of fractions of Poisson prime
factors of the algebras $k_\bfla[x_1,\dots,x_n]$. The method is a Poisson version
of the proof of \cite{\GLet, Theorem 2.1}.

\proclaim{3.4\. Theorem} Let $A= k_\bfla[x_1,\dots,x_n]$ for some antisymmetric matrix $\bfla \in
M_n(k)$, and let $P$ be a Poisson prime ideal of $A$. Then there exist a field extension
$K\supseteq k$ and an antisymmetric matrix $\bfmu\in M_m(k)$, for some $m\le n$, such that $\Fract
A/P \cong K_\bfmu(y_1,\dots,y_m)$ \pup{as Poisson algebras}.

In fact, $\bfmu$ is the upper left $m\times m$ submatrix of
$\bfsig\bfla\bfsig^{\tr}$, for some $\bfsig\in GL_n(\ZZ)$.
\endproclaim

\demo{Proof} Write $\bfla= (\lambda_{ij})$.

If $I= \{i\in \{1,\dots,n\} \mid x_i\notin P\}$ and $\bfla'$ is the submatrix of
$\bfla$ consisting of the rows and columns indexed by $I$, then there is a
Poisson prime ideal $P'$ in $k_{\bfla'}[x_i\mid i\in I]$ such that $A/P\cong
k_{\bfla'}[x_i\mid i\in I]/P'$. Thus, there is no loss of generality in assuming
that $x_i\notin P$ for all $i$. Now $P$ induces a Poisson prime ideal $Q$ in the
algebra $B= k_\bfla[x_1^{\pm1},\dots,x_n^{\pm1}]$ such that $Q\cap A= P$, and
$\Fract A/P\cong \Fract B/Q$ (as Poisson algebras).

Write $B= k_b\Gamma$ as in \S3.1, where $\Gamma= \ZZ^n$ and $b$ is the $k$-valued
antisymmetric bilinear form on $\Gamma$ obtained from $\bfla$. Set
$$\Gamma_Q= \{\alpha \in\Gamma \mid x^\alpha+Q \in Z_p(B/Q) \},$$
and observe that $\Gamma_Q$ is a subgroup of $\Gamma$. We claim that
$\Gamma/\Gamma_Q$ is torsionfree. If $\alpha\in \Gamma$ and $t\alpha\in \Gamma_Q$
for some $t\in \NN$, then $(x^\alpha+Q)^t= x^{t\alpha}+Q$ lies in $Z_p(B/Q)$,
whence
$$t(x^\alpha+Q)^{t-1} \{x^\alpha+Q,\, B/Q\}= \{(x^\alpha+Q)^t,\, B/Q\}= 0.$$
Since $t(x^\alpha+Q)^{t-1}$ is a unit in $B/Q$, it follows that $\{x^\alpha+Q,\,
B/Q\}= 0$, that is, $\alpha\in \Gamma_Q$. Therefore $\Gamma/\Gamma_Q$ is
torsionfree, as claimed.

Now $\Gamma/\Gamma_Q$ is a free abelian group of rank $m\le n$, so there exists a basis
$(\eps_1,\dots,\eps_n)$ for $\Gamma$ such that $\Gamma_Q$ is generated by
$\{\eps_{m+1},\dots,\eps_n\}$. Let $(\gamma_1,\dots,\gamma_n)$ denote the standard basis for
$\Gamma$. There is a matrix $\bfsig= (\sigma_{rs}) \in GL_n(\ZZ)$ such that each $\eps_r=
\sum_{s=1}^n \sigma_{rs} \gamma_s$, and we set $\bfxi$ equal to the antisymmetric matrix
$\bfsig\bfla\bfsig^{\tr} \in M_n(k)$. Then $B= k_\bfxi[z_1^{\pm1}, \dots, z_n^{\pm1}]$, where
$z_i= x^{\eps_i}$ for $i=1,\dots,n$. Hence, after replacing the $x_i$ and $\bfla$ by the $z_i$ and
$\bfxi$, we may assume that $(\eps_1,\dots,\eps_n)$ is the standard basis for $\Gamma$. In
particular, this means that $x_i+Q \in Z_p(B/Q)$ for $i=m+1,\dots,n$.

Set $K= \Fract Z_p(B/Q) \subseteq \Fract B/Q$, let $\bfmu\in M_m(k)$ be the
(antisymmetric) upper left $m\times m$ submatrix of $\bfla$, and
form the Poisson
$K$-algebra $C= K_\bfmu[y_1^{\pm1}, \dots, y_m^{\pm1}]$. There is a $K$-algebra
homomorphism $\phi:C \rightarrow \Fract B/Q$ such that $\phi(y_i)= x_i+Q$ for
$i=1,\dots,m$, and $\phi$ is a Poisson homomorphism because
$$\{\phi(y_i), \phi(y_j)\}= \{x_i,x_j\}+Q= \lambda_{ij}x_ix_j+Q= \phi(
\{y_i,y_j\})$$
for $i,j=1,\dots,m$. Since $x_i+Q\in K$ for $i>m$, the image of $\phi$ contains
all the $x_i+Q$, and so
$$\Fract \phi(C)= \Fract B/Q \cong \Fract A/P.$$
Thus, it only remains to show that $\phi$ is injective.

We claim that $C$ is Poisson simple. Identify $C$ with $K_c\Delta$ in the
notation of \S3.1, where $\Delta= \ZZ^m$ and $c$ is the $k$-valued antisymmetric
bilinear form on $\Delta$ obtained from $\bfmu$. Further, identify $\Delta$ with
the subgroup of $\Gamma$ generated by $\eps_1,\dots,\eps_m$; then $\Gamma= \Delta
\oplus \Gamma_Q$ and
$c$ is the restriction of $b$ to $\Delta\times\Delta$. We use Lemma 3.2 to prove
that
$Z_p(C)=K$, after which Corollary 3.3 will imply that $C$ is Poisson simple.
Thus, let $\alpha\in \Delta_c$, that is, $\alpha\in \Delta$ and $c(\alpha,-)
\equiv 0$. For $j=1,\dots,m$, we obtain
$$0= c(\alpha,\eps_j)= \sum_{i=1}^m \alpha_i\lambda_{ij},$$
and consequently
$$\{x^\alpha,x_j\}= \sum_{i=1}^m \alpha_i\lambda_{ij} x^\alpha x_j = 0.$$
Since $\{x^\alpha,x_j\} \in Q$ for $j>m$ (because $x_j+Q \in Z_p(B/Q)$), it follows that $x^\alpha+Q \in
Z_p(B/Q)$, and so $\alpha\in \Gamma_Q$. However, $\Delta \cap\Gamma_Q =0$,
forcing $\alpha=0$. We have proved that $\Delta_c =0$, and hence $Z_p(C)=K$ by
Lemma 3.2. Corollary 3.3 now implies that $C$ is Poisson simple, as claimed.

Since $\ker\phi$ is a Poisson ideal of $C$, it must be zero. Therefore $\phi$ is
injective, and the proof is complete. \qed\enddemo

We next construct a Poisson version of the derivation-deleting map introduced by
Cauchon in \cite{\Cau, Section 2}.

\proclaim{3.5\. Lemma} Let $A= B[x;\alpha,\delta]_p$ be a Poisson
polynomial algebra. Assume that $\delta$ is locally nilpotent, and that
$\alpha\delta= \delta(\alpha+s)$ for some $s\in\kx$. Then the rule
$$\theta(b)= \sum_{n=0}^\infty \dfrac1{n!}\left( \dfrac{-1}s \right)^n
\delta^n(b) x^{-n} \tag3-1$$
defines a $k$-algebra homomorphism $\theta: B\rightarrow
B[x^{\pm1}]$, and
$$\{x,\theta(b)\}= \theta\alpha(b)x \tag3-2$$
for all $b\in B$. \endproclaim

\demo{Proof} Note first that (3-1) at least defines a $k$-linear map $\theta:
B\rightarrow B[x^{\pm1}]$, and that $\theta(1)=1$. We compute that
$$\align \theta(ab) &= \sum_{n=0}^\infty \dfrac1{n!}\left( \dfrac{-1}s \right)^n
\delta^n(ab) x^{-n} = \sum_{n=0}^\infty \dfrac1{n!}\left( \dfrac{-1}s \right)^n
\sum_{l=0}^n {n\choose l} \delta^l(a) \delta^{n-l}(b) x^{-n} \\
&= \sum_{l,m=0}^\infty \dfrac1{l!} \dfrac1{m!} \left( \dfrac{-1}s \right)^{l+m}
\delta^l(a) \delta^m(b) x^{-l-m} =\theta(a)\theta(b)
\endalign$$
for all $a,b\in A$. Therefore $\theta$ is a $k$-algebra homomorphism.

For $b\in B$, we have
$$\align \{x,\theta(b)\} &= \sum_{n=0}^\infty \dfrac1{n!}\left( \dfrac{-1}s
\right)^n \bigl( \alpha\delta^n(b)x+ \delta^{n+1}(b) \bigr) x^{-n} \\
 &= \theta\delta(b)+ \sum_{n=0}^\infty \dfrac1{n!}\left( \dfrac{-1}s
\right)^n \delta^n(\alpha+ns)(b) x^{1-n} \\
 &= \theta\delta(b)+ \theta\alpha(b)x- \sum_{n=1}^\infty \dfrac1{(n-1)!} \left(
\dfrac{-1}s \right)^{n-1} \delta^n(b)x^{1-n} \\
 &= \theta\delta(b)+ \theta\alpha(b)x- \theta\delta(b)=
\theta\alpha(b)x. \endalign$$
This proves (3-2).  \qed\enddemo

\proclaim{3.6\. Lemma} Let $A= B[x;\alpha,\delta]_p$ be a Poisson
polynomial algebra, and assume that $\alpha\delta= \delta(\alpha+s)$ for some
$s\in k$. Then
$$\delta^n \bigl( \{a,b\} \bigr)= \sum_{l+m=n} {n\choose l} \bigl(
\{\delta^l(a),\delta^m(b)\}+ m\delta^l\alpha(a)\delta^m(b)-
l\delta^l(a)\delta^m\alpha(b) \bigr)  \tag3-3$$
for all $a,b\in B$ and $n\ge 0$. \endproclaim

\demo{Proof} Let $L$ denote the $k[x]$-linear map $\{x,-\} :A \rightarrow A$. Because of the
Jacobi identity for the Poisson bracket, $L$ is a Poisson derivation on $A$. Hence, $L$ satisfies
the Leibniz Rule
$$L^n\bigl( \{a,b\} \bigr)= \sum_{l+m=n} {n\choose l} \{L^l(a),L^m(b)\}  \tag3-4$$
for $a,b\in A$ and $n\ge 0$. Next, write $\equiv_1$ and $\equiv_2$ for congruence modulo the
ideals $(x)$ and $(x^2)$ in $A$, respectively. We claim that
$$L^n(a) \equiv_2 \delta^n(a)+ \bigl[ n\delta^{n-1}\alpha(a)+ \dfrac{n(n-1)s}{2}
\delta^{n-1}(a) \bigr] x  \tag3-5$$
for $a\in B$ and $n\ge 0$. This holds trivially when $n=0$,
and by construction when $n=1$.

If (3-5) holds for some $a\in B$ and some $n$, then
$$\align L^{n+1}(a) &\equiv_2 \delta^{n+1}(a)+ \alpha\delta^n(a)x+  \bigl[ n\delta^n\alpha(a)+
\dfrac{n(n-1)s}{2} \delta^n(a) \bigr] x \\
 &\equiv_2 \delta^{n+1}(a)+ \bigl[ (n+1)\delta^n\alpha(a)+ \dfrac{n(n+1)s}{2}
\delta^n(a) \bigr] x,  \endalign$$ because $\alpha\delta^n= \delta^n(\alpha+ns)$. Thus, by
induction, (3-5) holds.

Now let $a,b\in B$ and $n\ge 0$. In view of (3-4) and (3-5), we have
$$\alignedat1 L^n\bigl( \{a,b\} \bigr) &\equiv_1 \sum_{l+m=n} {n\choose l} \bigl\{
\delta^l(a)+ \bigl[ l\delta^{l-1}\alpha(a)+ \dfrac{l(l-1)s}{2}
\delta^{l-1}(a) \bigr] x,  \\
 &\qquad\qquad\qquad\qquad\qquad \delta^m(b)+ \bigl[ m\delta^{m-1}\alpha(b)+
\dfrac{m(m-1)s}{2} \delta^{m-1}(b) \bigr] x \bigr\} \\
 &\equiv_1  \sum_{l+m=n} {n\choose l} \bigl( \{\delta^l(a),\delta^m(b)\}+ U_{lm}+
V_{lm} \bigr),  \endalignedat  \tag3-6$$
 where
$$\align U_{lm} &= l\delta^{l-1}\alpha(a)\delta^{m+1}(b)-
m\delta^{m-1}\alpha(b)\delta^{l+1}(a)  \\
V_{lm} &= \dfrac{l(l-1)s}{2}
\delta^{l-1}(a)\delta^{m+1}(b)- \dfrac{m(m-1)s}{2}
\delta^{m-1}(b)\delta^{l+1}(a).  \endalign$$
Observe that
$$\alignedat1  \sum_{l+m=n} {n\choose l} U_{lm} &= \sum \Sb l+m=n\\ m>0 \endSb
{n\choose{l+1}}(l+1) \delta^l\alpha(a)\delta^m(b)  \\
 &\qquad -\sum \Sb l+m=n\\ l>0 \endSb
{n\choose{l-1}}(m+1) \delta^l(a)\delta^m\alpha(b)  \\
 &= \sum_{l+m=n} {n\choose l} \bigl( m\delta^l\alpha(a)\delta^m(b)-
l\delta^l(a)\delta^m\alpha(b) \bigr), \endalignedat  \tag3-7$$
while
$$\alignedat1  \sum_{l+m=n} {n\choose l} V_{lm} &= \dfrac{s}{2} \sum \Sb l+m=n\\ m>0 \endSb
{n\choose{l+1}}(l+1)l \delta^l(a)\delta^m(b)  \\
 &\qquad - \dfrac{s}{2} \sum \Sb l+m=n\\ l>0 \endSb
{n\choose{l-1}}(m+1)m \delta^l(a)\delta^m(b) \\
 &=0. \endalignedat  \tag3-8$$
 Combining (3-6), (3-7),
and (3-8), we obtain
$$ L^n\bigl( \{a,b\} \bigr)  \equiv_1
 \sum_{l+m=n} {n\choose l} \bigl(
 \{\delta^l(a),\delta^m(b)\}+ m\delta^l\alpha(a)\delta^m(b)-
l\delta^l(a)\delta^m\alpha(b) \bigr).  \tag3-9$$
Since $L^n\bigl( \{a,b\} \bigr) \equiv_1
\delta^n\bigl( \{a,b\} \bigr) $, (3-9) implies (3-3). \qed\enddemo

\proclaim{3.7\. Lemma} Under the hypotheses of Lemma {\rm 3.5}, the map $\theta$
is a Poisson homomorphism from $B$ to
$B[x^{\pm1};\alpha,\delta]_p$.
\endproclaim

\demo{Proof} For $a,b\in B$, we compute, using (3-2), that
$$\aligned  \{\theta(a),\theta(b)\} &= \sum_{l=0}^\infty \dfrac1{l!} \left(
\dfrac{-1}s \right)^l \{\delta^l(a)x^{-l}, \theta(b)\}  \\
 &= \sum_{l=0}^\infty \dfrac1{l!} \left(
\dfrac{-1}s \right)^l \bigl( \{\delta^l(a),\theta(b)\}x^{-l} - l
\delta^l(a)\{x,\theta(b)\}x^{-l-1} \bigr)  \\
 &= \sum_{l=0}^\infty \dfrac1{l!} \left(
\dfrac{-1}s \right)^l \bigl( \{\delta^l(a),\theta(b)\} -
l\delta^l(a)\theta\alpha(b) \bigr) x^{-l}  \\
 &= \sum_{l,m=0}^\infty \dfrac1{l!m!} \left(
\dfrac{-1}s \right)^{l+m} \bigl( \{\delta^l(a),\delta^m(b)x^{-m}\} -
l\delta^l(a)\delta^m\alpha(b)x^{-m} \bigr) x^{-l}  \\
  &= \sum_{l,m=0}^\infty \dfrac1{l!m!} \left(
\dfrac{-1}s \right)^{l+m} (C_{lm}+D_{lm}) x^{-l-m},
\endaligned  \tag3-10$$
where
$$\aligned C_{lm} & = \{\delta^l(a),\delta^m(b)\} +
m\delta^l\alpha(a)\delta^m(b)  - l\delta^l(a)\delta^m\alpha(b)  \\
D_{lm} &= lms \delta^l(a)\delta^m(b) + m \delta^{l+1}(a)\delta^m(b)x^{-1}
\endaligned  \tag3-11$$
for all $l$, $m$. The contribution of the $D_{lm}$ terms to the sum in (3-10) is
$$\aligned \sum_{l,m=0}^\infty \dfrac1{l!m!} \left(
\dfrac{-1}s \right)^{l+m} &\bigl[ lms \delta^l(a)\delta^m(b) x^{-l-m} +
m  \delta^{l+1}(a)\delta^m(b)x^{-l-m-1} \bigr]   \\
 &= \sum_{l,m=1}^\infty \dfrac{1}{(l-1)!(m-1)!} \dfrac{(-1)^{l+m}}{s^{l+m-1}}
\delta^l(a)\delta^m(b) x^{-l-m}  \\
 &\qquad + \sum_{t,m=1}^\infty \dfrac{1}{(t-1)!(m-1)!} \left(
\dfrac{-1}s \right)^{t-1+m} \delta^t(a)\delta^m(b) x^{-t-m}  \\
 &= 0.
\endaligned  \tag3-12$$
Because of (3-11), (3-12), and Lemma 3.6, we may simplify (3-10) to
$$\aligned  \{\theta(a),\theta(b)\} &= \sum_{l,m=0}^\infty \dfrac1{l!m!} \left(
\dfrac{-1}s \right)^{l+m} C_{lm} x^{-l-m}  \\
 &= \sum_{n=0}^\infty \dfrac1{n!} \left(
\dfrac{-1}s \right)^n \delta^n(\{a,b\}) x^{-n} = \theta(\{a,b\}).
\endaligned  $$
Therefore $\theta$ preserves the Poisson bracket. \qed
\enddemo

\proclaim{3.8\. Theorem} Under the hypotheses of Lemma {\rm 3.5}, the map
$\theta$  extends uniquely to an isomorphism of
Poisson Laurent polynomial algebras,
$$\theta: B[y^{\pm1};\alpha]_p @>{\;\cong\;}>> B[x^{\pm1};\alpha,\delta]_p\,,$$
such that $\theta(y)=x$. \endproclaim

\demo{Proof} First, $\theta$ extends uniquely to a $k$-algebra homomorphism
$B[y^{\pm1}] \rightarrow B[x^{\pm1}]$ such that $\theta(y)=x$. In view of (3-2)
and Lemma 3.7, the forms $\theta(\{-,-\})$ and $\{\theta(-),\theta(-)\}$ agree on
pairs of elements from $B\cup\{y^{\pm1}\}$, from which we see that they agree on
pairs of elements from $A$. Hence, the extended map $\theta$ is a Poisson
homomorphism, and it only remains to show that
$\theta$ is bijective.

For surjectivity, we already have $x^{\pm1}= \theta(y^{\pm1})$, and so we just
need to see that $B$ is contained in the image of $\theta$. Given $b\in B$,
there is some $l\ge0$ such that $\delta^l(b)=0$, and we proceed by induction on
$l$. If $l\le 1$, then $\delta(b)=0$ and
$\theta(b)=b$. Now let $l>1$, and write $\theta(b)= b+ \sum_{n=1}^{l-1} \lambda_n
\delta^n(b)x^{-n}$ for some $\lambda_n\in k$. Since $\delta^{l-1}(\delta^n(b))=0$
for $n=1,\dots,l-1$, we can assume by induction that
$\delta^1(b),\dots,\delta^{l-1}(b)$ are in the image of $\theta$. Consequently,
$\theta(b)-b$ is in the image of $\theta$, and thus $b$ is in the image of
$\theta$. This establishes the induction step, and proves that $\theta$ is
surjective.

Let $p\in B[y^{\pm1}]$ be nonzero, and write $p= \sum_{i=l}^m b_iy^i$ for some
$b_i\in B$ and some integers $l\le m$, with
$b_m\ne 0$. Each of the terms
$\theta(b_iy^i)$ is a Laurent polynomial of the form $b_ix^i + [\text{lower
terms}]$. Hence,
$\theta(p)= b_mx^m  + [\text{lower
terms}]$, and thus $\theta(p) \ne 0$. Therefore $\theta$ is injective.
\qed\enddemo

The following is the main result addressing our quadratic version of the Poisson Gelfand-Kirillov
problem. It is a Poisson version of \cite{\Cau, Th\'eor\`eme 6.1.1}.

\proclaim{3.9\. Theorem} Let $A= k[x_1] [x_2;\alpha_2,\delta_2]_p \cdots
[x_n;\alpha_n,\delta_n]_p$ be an iterated Poisson polynomial algebra such that

{\rm (a)} $\delta_i$ is locally nilpotent for all $i$.

{\rm (b)} There exist $s_i\in\kx$ such that $\alpha_i\delta_i=
\delta_i(\alpha_i+s_i)$ for all $i$.

{\rm (c)} There exist $\lambda_{ij}\in k$ such that $\alpha_i(x_j)=
\lambda_{ij}x_j$ for all $i>j$.

\noindent Let $\bfla= (\lambda_{ij})$ be the antisymmetric matrix in $M_n(k)$ whose entries below
the diagonal agree with the scalars in {\rm (c)}. Then:

{\rm (1)} $\Fract A \cong k_{\bfla}(y_1,\dots,y_n)$.

{\rm (2)} For any Poisson prime ideal $P$ of $A$, there exist a field extension $K\supseteq k$ and
an antisymmetric matrix $\bfmu\in M_m(k)$, for some $m\le n$, such that $\Fract A/P \cong
K_\bfmu(y_1,\dots,y_m)$ \pup{as Poisson algebras}. In fact, $\bfmu$ is the upper left $m\times m$
submatrix of $\bfsig\bfla\bfsig^{\tr}$, for some $\bfsig\in GL_n(\ZZ)$.
\endproclaim

\demo{Proof} Let $P$ be an arbitrary Poisson prime ideal of $A$, and set $B=
k_\bfla[z_1,\dots,z_n]$. In view of Theorem 3.4, it suffices to show that \roster
\item"(*)" $\Fract A/P \cong \Fract B/Q$ for some Poisson prime ideal $Q$ of $B$, where $Q=0$ if
$P=0$.
\endroster
We prove (*) via a triple induction: first, with respect to $n$; second, with respect to the
number $d$ of indices $i$ for which $\delta_i\ne 0$; and third (downward), with respect to the
maximum index $t$ for which $\delta_t \ne 0$. (If $d=0$, we take $t=n+1$.) Since there is nothing
to prove if $n=1$ or $t=n+1$, we may assume that $n\ge 2$ and $t\le n$.

{\bf Case 1}: $x_n\in P$. Then there exists a Poisson prime ideal $P'$ in the
algebra
$$A'= k[x_1] [x_2;\alpha_2,\delta_2]_p \cdots
[x_{n-1};\alpha_{n-1},\delta_{n-1}]_p$$ such that $A/P \cong A'/P'$. By our primary induction,
there is a Poisson prime ideal $Q'$ in the algebra $B'= k_{\bfla'}[z_1,\dots,z_{n-1}]$, where
$\bfla'$ is the upper left $(n-1)\times(n-1)$ submatrix of $\bfla$, such that $\Fract A'/P' \cong
\Fract B'/Q'$. Observe that $Q= Q'+Bz_n$ is a Poisson prime ideal of $B$ such that $B'/Q' \cong
B/Q$. Thus,  $\Fract A/P \cong \Fract B/Q$.

{\bf Case 2}: $x_n\notin P$ and $t=n$. Then $\delta_n \ne 0$. Set
$$A'= k[x_1] [x_2;\alpha_2,\delta_2]_p \cdots
[x_{n-1};\alpha_{n-1},\delta_{n-1}]_p[y;\alpha_n]_p.$$ By Theorem 3.8, $A[x_n^{-1}] \cong
A'[y^{-1}]$, and so there exists a Poisson prime ideal $P'$ in $A'$ such that $\Fract A/P \cong
\Fract A'/P'$, where $P'=0$ if $P=0$. The number of nonzero maps among
$\delta_2,\dots,\delta_{n-1}$ is $d-1$. Thus, our secondary induction yields (*) in this case.

{\bf Case 3}: $t<n$. Then $\delta_n =0$. Since $\{x_n,x_1\}= \lambda_{n1}x_1x_n$, we see that
$\{x_n,k[x_1]\} \subseteq k[x_1]x_n$, and so $k[x_1,x_n]$ is a Poisson polynomial algebra of the
form $k[x_1][x_n;\alpha'_n]_p$. For $i=2,\dots,n-1$, we have
$$\{x_i,k[x_1,\dots,x_{i-1}]\} \subseteq k[x_1,\dots,x_{i-1}]x_i+
k[x_1,\dots,x_{i-1}]$$
and $\{x_i,x_n\}= -\lambda_{ni}x_ix_n= \lambda_{in}x_nx_i$, from which it
follows that
$$\{x_i,k[x_1,\dots,x_{i-1},x_n]\} \subseteq k[x_1,\dots,x_{i-1},x_n]x_i+
k[x_1,\dots,x_{i-1},x_n].$$
Hence, we may rewrite $A$ in the form
$$A= k[x_1][x_n;\alpha'_n]_p [x_2;\alpha'_2,\delta'_2]_p \cdots
[x_{n-1};\alpha'_{n-1},\delta'_{n-1}]_p$$
for suitable $\alpha'_i$ and $\delta'_i$, such that
$\alpha'_i(x_j)= \lambda_{ij}x_j$ for $j<i$ and for $j=n$. Note that $\alpha'_i$ and $\delta'_i$
restrict to $\alpha_i$ and $\delta_i$ on $k[x_1,\dots,x_{i-1}]$, and that $\delta'_i(x_n)=0$. It
follows easily that $\delta'_i$ is locally nilpotent, and that $\alpha'_i\delta'_i=
\delta'_i(\alpha'_i+s_i)$. Finally, the map $\delta'_t$ is nonzero because it restricts to
$\delta_t$, and it occurs in position $t+1$ in the list $0,0,\delta'_2,\dots,\delta'_{n-1}$. Thus,
our tertiary induction yields (*) in this case.

Therefore (*) holds in all cases, and the theorem is proved. \qed\enddemo

\head 4. Poisson polynomial algebras satisfying the quadratic Gelfand-Kirillov property \endhead

We apply Theorem 3.9 to the algebras discussed in Section 2, to obtain the following result.

\proclaim{4.1\. Theorem} Let $A$ be any of the Poisson algebras of \S\S{\rm 2.2(b), 2.3(b),
2.4(b), 2.5(b), 2.6(a), 2.7(a)}. For any Poisson prime ideal $P$ of $A$, there exist a field
extension $K\supseteq k$ and an antisymmetric matrix $\bfmu\in M_m(k)$, for some $m$, such that
$\Fract A/P \cong K_\bfmu(y_1,\dots,y_m)$ \pup{as Poisson algebras}. In case $P=0$, we have $K=k$
and $m= \trdeg_k A$. \qed\endproclaim

In each case, bounds on $m$ and restrictions on $\bfmu$ can be obtained via Theorem 3.9. We leave
details to the interested reader. The form of the Poisson field $\Fract A$ is given below.

The examples in Section 2 are already expressed as iterated Poisson polynomial algebras, and so
what remains is to establish hypotheses (a), (b), (c) of Theorem 3.9 in each case. For (a), the
following observation is helpful: To check that a derivation $\delta$ on an algebra $A$ is locally
nilpotent, it suffices to check that $\delta$ is locally nilpotent on a set of algebra generators
for $A$. (This follows directly from the Leibniz Rule for $\delta$.) Hypothesis (b) is built into
the situation of Corollary 1.8, as follows.

\proclaim{4.2\. Lemma} Let $A= k[x_1] [x_2;\alpha_2,\delta_2]_p \cdots [x_n;\alpha_n,\delta_n]_p$
be an iterated Poisson polynomial algebra, supporting a rational Poisson action by a torus $H$
such that $x_1,\dots,x_n$ are $H$-eigen\-vec\-tors. Assume that there exist
$\eta_1,\dots,\eta_n\in \hfrak= \Lie H$ such that $\eta_i.x_j= \alpha_i(x_j)$ for $i>j$ and the
$\eta_i$-eigenvalue of $x_i$, call it $s_i$, is nonzero for each $i$. Then $\alpha_i\delta_i=
\delta_i(\alpha_i+s_i)$ for all $i$.  \endproclaim

\demo{Proof}  Fix $i\in \{2,\dots,n\}$. Since the derivations $\eta_i.(-)$ and $\alpha_i$ agree on
$x_1,\dots,x_{i-1}$, they must agree on the algebra $A_{i-1}= k[x_1,\dots,x_{i-1}]$. Let $y_i$
denote the $H$-eigenvalue of $x_i$, so that $x_i\in A_{y_i}$. Then $s_ix_i= \eta_i.x_i=
(\eta_i|y_i)x_i$, and so $(\eta_i|y_i)= s_i$.

Consider an $H$-eigenvector $f\in A_{i-1}$, say $f\in A_z$ for some $z\in X(H)$, and note that
$$\{x_i,f\}= \alpha_i(f)x_i+ \delta_i(f)= (\eta_i.f)x_i+ \delta_i(f)= (\eta_i|z)fx_i+
\delta_i(f).$$
As shown in the proof of Lemma 1.6, $\{x_i,f\}\in A_{y_i+z}$. Since also $fx_i\in
A_{y_i+z}$, we see that $\delta_i(f)\in A_{y_i+z}$. Consequently,
$$\align \alpha_i\delta_i(f) &= \eta_i.\delta_i(f)= (\eta_i|y_i+z)\delta_i(f)= \delta_i \bigl(
(\eta_i|z)f+
(\eta_i|y_i)f \bigr)  \\
 &= \delta_i(\eta_i.f+ s_if)= \delta_i(\alpha_i+s_i)(f).  \endalign$$
 The lemma then follows from the rationality of the action of the torus $H$.
 \qed\enddemo

 Since we have shown that the examples in Section 2 satisfy the hypotheses of Corollary 1.8, we
 conclude from Lemma 4.2 that they also satisfy hypothesis (b) of Theorem 3.9.

 \definition{4.3} The algebra $A$ of \S2.2(b) is just $k_{\bfq}[x_1,\dots,x_n]$, and Theorem 3.4
 applies.
 \enddefinition

 \definition{4.4} Let $A= \OMn$ with the Poisson bracket given in (2-2). As the case $\lambda=0$
 is covered by \S4.3, we assume that $\lambda \ne 0$. Condition (c) of Theorem 3.9 is given by
 (2-3). The maps $\delta_{lm}$ in this algebra satisfy
 $$\delta_{lm}(X_{ij})= \cases \lambda X_{im}X_{lj} &\qquad (l>i,\, m>j) \\  0 &\qquad
 \text{(otherwise)}. \endcases$$
 In particular, $\delta^2_{lm}(X_{ij}) =0$ for all $(i,j) \lexle (l,m)$, whence $\delta_{lm}$ is
 locally nilpotent. Thus, the hypotheses of Theorem 3.9 are satisfied. In particular, the theorem
 implies that $\Fract A \cong k(Y_{ij} \mid i,j=1,\dots,n)$ with
 $$\{ Y_{lm}, Y_{ij} \}= \cases (p_{li}+p_{jm})Y_{ij}Y_{lm} &\quad (l\ge i,\ m>j)\\
 (\lambda+ p_{li}+p_{jm})Y_{ij}Y_{lm} &\quad (l>i,\ m\le j).\endcases$$
 \enddefinition

 \definition{4.5} Let $A= A^{P,Q}_{n,\Gamma}(k)$ as in \S2.4(b). Condition (c) of Theorem 3.9 is
 given by (2-6). The maps $\delta_i$ here satisfy
 $$\align  \delta_i(x_j) &= \delta_i(y_j) =0 \qquad (j<i)  \\
 \delta_i(x_i) &= -\sum_{l<i} (q_l-p_l) x_ly_l \,.  \endalign$$
 Thus $\delta_i^2$ vanishes on $x_1,y_1,\dots, x_{i-1},y_{i-1}, x_i$, whence $\delta_i$ is locally
 nilpotent. In this case, Theorem 3.9 shows that $\Fract A \cong k(v_1,w_1,\dots,v_n,w_n)$ with
$$\alignedat2 \{w_i,w_j\} & = \gamma_{ij}w_iw_j &\qquad\qquad\qquad&(\text{all\
} i,j) \\
\{v_i,w_j\} &= (p_j+\gamma_{ji})v_iw_j &&(i<j) \\
\{v_i,w_j\} &= (q_j+ \gamma_{ji})v_iw_j && (i\ge j) \\
\{v_i,v_j\} &= (q_i-p_j+\gamma_{ij})v_iv_j && (i<j).\endalignedat$$
 \enddefinition

 \definition{4.6} Let $A$ be as in \S2.5(b). Condition (c) of Theorem 3.9 is clear from (2-8). The
 maps $\delta_i$ here satisfy
 $$\align  \delta_i(x_j) &= \delta_i(y_j)= \delta_i(z_0) =0 \qquad (j<i)  \\
 \delta_i(x_i) &= -\sum_{l<i} (q_l-p_l) x_ly_l- \lambda z_0^2 \,.  \endalign$$
 Thus $\delta_i^2$ vanishes on $z_0,x_1,y_1,\dots, x_{i-1},y_{i-1}, x_i$, whence $\delta_i$ is locally
 nilpotent. We see from Theorem 3.9 that $\Fract A \cong k(u_0,v_1,w_1,\dots,v_n,w_n)$ with
$$\alignedat2 \{u_0,v_i\} &= -(p_i/2)u_0v_i  &&(\text{all\ } i)  \\
\{u_0,w_i\} &= (p_i/2)u_0w_i  &&(\text{all\ } i)  \\
\{w_i,w_j\} & = \gamma_{ij}w_iw_j &\qquad\qquad\qquad&(\text{all\
} i,j) \\
\{v_i,w_j\} &= (p_j+\gamma_{ji})v_iw_j &&(i<j) \\
\{v_i,w_j\} &= (q_j+ \gamma_{ji})v_iw_j && (i\ge j) \\
\{v_i,v_j\} &= (q_i-p_j+\gamma_{ij})v_iv_j && (i<j).\endalignedat$$
 \enddefinition

\definition{4.7} Let $A$ be as in \S2.6(a). Condition (c) of Theorem 3.9 is given by (2-10). The
maps $\delta_{lm}$ in this algebra satisfy
$$\delta_{lm}(y_{ij})= \cases -2y_{im}y_{lj} &\qquad (i<l\le j<m) \\  -2y_{il}y_{jm} -2y_{im}y_{jl}
&\qquad (i\le j<l\le m) \\  0 &\qquad \text{(otherwise)}  \endcases$$ for $l\le m$ and $i\le j$
with $(i,j)\lexle (l,m)$. It follows that $\delta^3_{lm}(y_{ij})=0$ for all $(i,j)\lexle (l,m)$,
whence $\delta_{lm}$ is locally nilpotent. In this case, Theorem 3.9 implies that $\Fract A \cong
k(z_{ij} \mid 1\le i\le j\le n)$ with
$$\{z_{ij},z_{lm}\}= \cases  z_{ij}z_{lm}  &\quad \bigl( (i=l<j<m)
\text{\ or\ } (i<l<j=m) \text{\ or\ } (i<j=l<m) \bigr)  \\
2z_{ij}z_{lm}  &\quad \bigl( (i=j=l<m) \text{\ or\ } (i<j=l=m) \bigr)  \\
0  &\quad (\text{otherwise})  \endcases$$
for $l\le m$ and $i\le j$ with $(i,j) \lexle (l,m)$.
\enddefinition

\definition{4.8} Let $A$ be as in \S2.7(a). Condition (c) of Theorem 3.9 is given by (2-12). The
maps $\delta_{lm}$ in this algebra satisfy
$$\delta_{lm}(y_{ij})= \cases -2y_{im}y_{lj} &\qquad (i<l< j<m) \\  -2y_{il}y_{jm} +2y_{im}y_{jl}
&\qquad (i< j<l< m) \\  0 &\qquad \text{(otherwise)}  \endcases$$ for $l< m$ and $i< j$ with
$(i,j)\lexle (l,m)$. It follows that $\delta^2_{lm}(y_{ij})=0$ for all $(i,j)\lexle (l,m)$, whence
$\delta_{lm}$ is locally nilpotent. In this case, finally, we see from Theorem 3.9 that $\Fract A
\cong k(z_{ij} \mid 1\le i< j\le n)$ with
$$\{z_{ij},z_{lm}\}= \cases  z_{ij}z_{lm}  &\qquad (\text{if\ } {\bigm|} \{i,j\}\cap
\{l,m\} {\bigm|}= 1)  \\   0  &\qquad (\text{otherwise})  \endcases$$
for $l<m$ and $i<j$ with
$(i,j) \lexle (l,m)$.
\enddefinition

\head 5. Isomorphism invariants of quadratic Poisson fields \endhead

In this final section of the paper, we address the question of when Poisson fields
$k_{\bfla}(x_1,\dots,x_n)$ and $k_{\bfmu}(x_1,\dots,x_n)$ can be isomorphic. It is easily seen
that a sufficient condition is the existence of an invertible integer matrix $A$ such that $\bfmu=
A\bfla A^{\tr}$ (Lemma 5.1), and we show that in a number of cases, this condition is also necessary. The
method is to show that the set of matrices $B\bfla B^{\tr}$, for $B\in M_n(\ZZ)$, is an invariant
of $k_{\bfla}(x_1,\dots,x_n)$. By similar means, we also show that $k_{\bfla}(x_1,\dots,x_n)$
cannot be isomorphic to any Poisson-Weyl field. The invariants we use are Poisson analogs of some
invariants introduced by Alev and Dumas in \cite{\AlDu}.

For purposes of computation in $k_{\bfla}(x_1,\dots,x_n)$, observe that the
Poisson bracket of any monomials $x^a$ and $x^b$ is given by
$$\{x^a,x^b\}_{\bfla}= \sum_{l,m=1}^n a_lb_m\lambda_{lm}x^{a+b}= (a\bfla b^{\tr})
x^{a+b},  \tag5-1$$
where $a,b\in \ZZ^n$ are viewed as row vectors.

\proclaim{5.1\. Lemma} Let $\bfla,\bfmu\in M_n(k)$ be antisymmetric, and assume
there exists $A\in GL_n(\ZZ)$ such that $\bfmu= A\bfla A^{\tr}$. Then
$k_{\bfla}(x_1,\dots,x_n) \cong k_{\bfmu}(x_1,\dots,x_n)$ \pup{as Poisson algebras
over $k$}.  \endproclaim

\demo{Proof} Let $a_1,\dots,a_n$ denote the rows of $A$, set $y_i= x^{a_i}$ for
$i=1,\dots,n$, and observe using (5-1) that
$$\{y_i,y_j\}_{\bfla}= (a_i\bfla a_j^{\tr}) y_iy_j= \bfmu_{ij} y_iy_j  \tag5-2$$
for all $i$, $j$. Since $A$ is invertible, $x_1,\dots,x_n$ all lie in $k(y_1,\dots,y_n)$, so the
$y_i$ are algebraically independent over $k$ and $k(y_1,\dots,y_n)= k(x_1,\dots,x_n)$. Hence,
there is a $k$-algebra automorphism $\phi$ of $k(x_1,\dots,x_n)$ sending $y_i\mapsto x_i$ for all
$i$. Since the Poisson brackets $\{-,-\}_{\bfla}$ and $\{-,-\}_{\bfmu}$ are determined by the
values $\{y_i,y_j\}_{\bfla}$ and $\{x_i,x_j\}_{\bfmu}$, we conclude that $\phi$ is a Poisson
isomorphism of $k_{\bfla}(x_1,\dots,x_n)$ onto $k_{\bfmu}(x_1,\dots,x_n)$.  \qed\enddemo

\proclaim{5.2\. Proposition} Let $K=k_{\bfla}(x_1,\dots,x_n)$ for some
antisymmetric $\bfla\in M_n(k)$.

{\rm (a)} If $B_{\bfla}$ is the $k$-subspace of $K$ spanned by $\{\{f,g\} \mid f,g\in K\}$, then
$B_{\bfla}\cap k= \{0\}$.

{\rm (b)} For any $n$-tuple $y= (y_1,\dots,y_n)$ of nonzero elements of $K$, let $C_{\bfla}(y)$
denote the matrix $\bigl( \{y_i,y_j\}(y_iy_j)^{-1} \bigr) \in M_n(K)$. If $C_{\bfla}=
\{C_{\bfla}(y)\mid y\in (K^\times)^n\}$, then $C_{\bfla}\cap M_n(k)= \{A\bfla A^{\tr} \mid A\in
M_n(\ZZ)\}$.  \endproclaim

\demo{Proof} Put the lexicographic order on $\ZZ^n$, and let $L$ denote the corresponding
Hahn-Laurent power series field in $x_1,\dots,x_n$ (cf\. \cite{\Cohn, Theorem VII.3.8}; a more
detailed treatment can be found in \cite{\DiLe, Section 2}). The field $L$ consists of formal
series $\sum_{a\in I} \alpha_ax^a$ where $I$ is a well-ordered subset of $\ZZ^n$ and the
$\alpha_a\in k$. Finite sums in $L$ are identified with Laurent polynomials in
$k[x_1^{\pm1},\dots,x_n^{\pm1}]$. Since $L$ is a field, it thus contains (a copy of) $K$. Let
$\pi: L\rightarrow k$ be the $k$-linear map that gives the constant term (i.e., the coefficient of
$x^{\bfz}$) of elements of $L$. Observe that the Poisson bracket on $K$ extends to $L$ by setting
$$\{f,g\}= \sum_{i,j=1}^n \lambda_{ij}x_ix_j \dfrac{\partial f}{\partial
x_i}\dfrac {\partial g}{\partial x_j}$$
for $f,g\in L$. This formula gives a well-defined element
of $L$ because the supports of $x_i (\partial f/\partial x_i)$ and $x_j (\partial g/\partial x_j)$
are contained in those of $f$ and $g$.

(a) It suffices to show that $\pi(\{f,g\})=0$ for any $f,g\in L$. Write $f=
\sum_{a\in I} \alpha_ax^a$ and $g= \sum_{b\in J} \beta_bx^b$ where $I$ and $J$ are
well-ordered subsets of $\ZZ^n$ and the $\alpha_a,\beta_b\in k$. Then
$$\{f,g\}= \sum_{i,j=1}^n \lambda_{ij} \bigl( \sum_{a\in I} a_i\alpha_ax^a \bigr) \bigl(
\sum_{b\in J} b_j\beta_bx^b \bigr)= \sum_{a\in I,\, b\in J} \bigl( \sum_{i,j=1}^n
\lambda_{ij}a_ib_j \bigr) \alpha_a\beta_bx^{a+b},  \tag5-3$$
and consequently
$$\pi\bigl( \{f,g\} \bigr)=  \sum \Sb a\in I,\, b\in J\\ a+b=\bfz \endSb \bigl(
\sum_{i,j=1}^n\lambda_{ij}a_ib_j \bigr) \alpha_a \beta_b= \sum_{a\in I\cap(-J)}
\bigl( -\sum_{i,j=1}^n\lambda_{ij}a_ia_j \bigr) \alpha_a \beta_{-a}.$$
Since $\bfla$ is antisymmetric, each of the sums
$\sum_{i,j=1}^n\lambda_{ij}a_ia_j$ is zero, and thus $\pi(\{f,g\})=0$, as desired.

(b) It follows from (5-2) that $A\bfla A^{\tr} \in C_{\bfla}$ for all $A\in M_n(\ZZ)$. Hence, it
suffices to show that for any $(y_1,\dots,y_n) \in (L^\times)^n$, the matrix $\bigl(\pi
(\{y_i,y_j\}(y_iy_j)^{-1}) \bigr)$ has the form $A\bfla A^{\tr}$ for some $A\in M_n(\ZZ)$.

Write each $y_i= \sum_{a\in I(i)} \alpha_{ia}y^a$ where $I(i)$ is a well-ordered subset of $\ZZ^n$
with minimum element $m(i)$, the $\alpha_{ia} \in k$, and $\alpha_{i,m(i)} \ne 0$. Note that
$y_i^{-1}= \sum_{b\in J(i)} \beta_{ib}x^b$ where $J(i)$ is a well-ordered subset of $\ZZ^n$ with
minimum element $-m(i)$, the $\beta_{ib}\in k$, and $\beta_{i,-m(i)}= \alpha_{i,m(i)}^{-1}$.

For any $i,j=1,\dots,n$, the series $\{y_i,y_j\}$ is supported on the set of those $c\in \ZZ^n$
for which $c\ge m(i)+m(j)$ (cf\. (5-3)), and so
$$\pi\bigl( \{y_i,y_j\}(y_iy_j)^{-1} \bigr)= \{x^{m(i)},x^{m(j)}\}x^{-m(i)-m(j)}=
m(i)\bfla m(j)^{\tr}$$
by (5-1). Thus, $\bigl( \pi( \{y_i,y_j\}(y_iy_j)^{-1} ) \bigr)= A\bfla
A^{\tr}$ where $A$ is the matrix in $M_n(\ZZ)$ with rows $m(1),\dots,\allowmathbreak m(n)$.
\qed\enddemo

The following corollaries give two immediate applications of Proposition 5.2. They are Poisson
analogs of results of Alev and Dumas, who proved that the quotient division ring of a quantum
plane $\O_q(k^2)$ cannot be isomorphic to a Weyl skew field \cite{\AlDu, Corollaire 3.11(a)}, and
that for nonroots of unity $q,r\in \kx$, the quotient division rings of $\O_q(k^2)$ and
$\O_r(k^2)$ are isomorphic if and only if $q= r^{\pm1}$ \cite{\AlDu, Corollaire 3.11(c)}.

\proclaim{5.3\. Corollary} Let $\bfla\in M_n(k)$ be antisymmetric. Then $k_{\bfla}(x_1,\dots,x_n)$
is not isomorphic to a Poisson-Weyl field. In fact, it is not isomorphic to any Poisson field
containing elements $x$ and $y$ with $\{x,y\}=1$.  \endproclaim

\demo{Proof} By Proposition 5.2(a), $\{x,y\} \ne 1$ for all $x,y\in
k_{\bfla}(x_1,\dots,x_n)$.  \qed\enddemo

\proclaim{5.4\. Corollary} Let $\bfla= \left[ \smallmatrix 0&\lambda\\ -\lambda&0
\endsmallmatrix \right]$ and $\bfmu= \left[ \smallmatrix 0&\mu\\ -\mu&0
\endsmallmatrix \right]$ for some $\lambda,\mu\in k$. Then $k_{\bfla}(x_1,x_2)
\cong k_{\bfmu}(x_1,x_2)$ if and only if $\lambda= \pm\mu$.  \endproclaim

\demo{Proof} If $\lambda=-\mu$, the $k$-algebra automorphism  of
$k(x_1,x_2)$ fixing
$x_1$ and sending $x_2\mapsto x_2^{-1}$ transforms $\{-,-\}_{\bfla}$ to
$\{-,-\}_{\bfmu}$, providing a Poisson
isomorphism of
$k_{\bfla}(x_1,x_2)$ onto $k_{\bfmu}(x_1,x_2)$.

Conversely, assume that $k_{\bfla}(x_1,x_2) \cong
k_{\bfmu}(x_1,x_2)$. By Proposition 5.2(b),
$$\{A\bfla A^{\tr} \mid A\in M_n(\ZZ)\}= \{B\bfmu B^{\tr} \mid B\in M_n(\ZZ)\},$$
from which we see that $\ZZ\lambda= \ZZ\mu$. Since $\chr k=0$, this implies
$\lambda= \pm\mu$.  \qed\enddemo

Cases (b) and (c) of the following theorem are Poisson analogs of results of Panov \cite{\Pan,
Theorem 2.19} and Richard \cite{\Ric, Th\'eor\`eme 4.2}.

\proclaim{5.5\. Theorem} Let $\bfla,\bfmu\in M_n(k)$ be antisymmetric, and assume
that one of the following holds:

{\rm (a)} $\bfla\in GL_n(k)$.

{\rm (b)} The subgroup $\sum_{i,j=1}^n \ZZ \lambda_{ij}$ of $(k,+)$ is cyclic.

{\rm (c)} The subgroup $\sum_{i,j=1}^n \ZZ \lambda_{ij}$ of $(k,+)$ is free
abelian of rank $n(n-1)/2$.

Then $k_{\bfla}(x_1,\dots,x_n) \cong k_{\bfmu}(x_1,\dots,x_n)$ if and
only if there exists $A\in GL_n(\ZZ)$ such that $\bfmu= A\bfla
A^{\tr}$. \endproclaim

\demo{Proof} Since the theorem is clear if $n=1$, we may assume that
$n\ge 2$. Sufficiency is given by Lemma 5.1. Conversely, assume that
$k_{\bfla}(x_1,\dots,x_n) \cong k_{\bfmu}(x_1,\dots,x_n)$. In view of Proposition
5.2(b), there exist $A,B\in M_n(\ZZ)$ such that $\bfmu= A\bfla
A^{\tr}$ and $\bfla= B\bfmu B^{\tr}$. Note that $\bfla= (BA)\bfla(BA)^{\tr}$.

(a) In this case, it follows from the equation $\bfla= (BA)\bfla(BA)^{\tr}$ that
$\det(BA)^2 =1$, and consequently $A,B\in GL_n(\ZZ)$.

(b) By assumption, $\sum_{i,j=1}^n \ZZ \lambda_{ij} =\ZZ\lambda$ for some
$\lambda\in k$. If $\lambda=0$, then $\bfla= \bfz$ and $\{-,-\}_{\bfla}$ vanishes.
In this case, $\{-,-\}_{\bfmu}$ must also vanish, whence $\bfmu=\bfz$ and $\bfmu=
I\bfla I^{\tr}$.

Now assume that $\lambda\ne 0$. Then $\lambda^{-1}\bfla$ is an
antisymmetric integer matrix, and so there exists $C\in GL_n(\ZZ)$ such that
$$C (\lambda^{-1}\bfla) C^{\tr}= \left[ \matrix 0&d_1 &0&\cdots &&&\cdots&0\\
-d_1&0 &&&&&&\vdots\\  0&&\ddots\\  \vdots&&&0&d_r\\
&&&-d_r&0\\  &&&&&0\\  \vdots&&&&&&\ddots&\vdots\\  0&\cdots&&&&&\cdots&0
\endmatrix \right]$$
for some nonzero integers $d_1,\dots,d_r$ (e.g., \cite{\New, Theorem IV.1}). Hence, we obtain a
block matrix decomposition
$$C \bfla C^{\tr}= \left[ \matrix \Lambda&\bfz\\ \bfz&\bfz \endmatrix \right]$$
with $\Lambda\in GL_{2r}(k)$. Since $C$ is invertible over $\ZZ$, we may replace $\bfla$ by $C
\bfla C^{\tr}$, and so there is no loss of generality in assuming that $\bfla= \left[ \matrix
\Lambda&\bfz\\ \bfz&\bfz \endmatrix \right]$.

The equations $\bfmu= A\bfla A^{\tr}$ and $\bfla= B\bfmu B^{\tr}$ imply that $\bfla$ and $\bfmu$
have the same rank, namely $2r$, and that $\sum_{i,j=1}^n \ZZ\mu_{ij}= \sum_{i,j=1}^n \ZZ
\lambda_{ij} =\ZZ\lambda$.
Hence, we also obtain a block matrix decomposition $D\bfmu D^{\tr}= \left[ \matrix M&\bfz\\
\bfz&\bfz \endmatrix \right]$ for some $D\in GL_n(\ZZ)$ and some $M\in GL_{2r}(k)$. As above,
there is no loss of generality in assuming that $\bfmu= \left[ \matrix M&\bfz\\ \bfz&\bfz
\endmatrix \right]$.

Write $A$ and $B$ in block form as
$$\xalignat2 A &= \left[ \matrix A_{11}&A_{12}\\ A_{21}&A_{22} \endmatrix \right]
&B &= \left[ \matrix B_{11}&B_{12}\\ B_{21}&B_{22} \endmatrix \right]
\endxalignat$$
where $A_{11}$ and $B_{11}$ are $2r\times 2r$. The equations $\bfmu= A\bfla
A^{\tr}$ and $\bfla= B\bfmu B^{\tr}$ now say that
$$\xalignat2 \left[ \matrix M&\bfz\\ \bfz&\bfz
\endmatrix \right] &= \left[ \matrix A_{11}LA_{11}^{\tr}
&A_{11}LA_{21}^{\tr} \\ A_{21}LA_{11}^{\tr} &A_{21}LA_{21}^{\tr} \endmatrix \right]  &\left[
\matrix \Lambda&\bfz\\ \bfz&\bfz
\endmatrix \right] &= \left[
\matrix B_{11}M B_{11}^{\tr} &B_{11}M B_{21}^{\tr} \\ B_{21}M B_{11}^{\tr} &B_{21}M B_{21}^{\tr}
\endmatrix \right],  \endxalignat$$
and so $\Lambda= (B_{11}A_{11})\Lambda (B_{11}A_{11})^{\tr}$. As in case (a), it follows that
$A_{11}\in GL_{2r}(\ZZ)$. Hence, the matrix $E= \left[ \matrix A_{11}&0\\
0&I_{n-2r} \endmatrix \right]$ lies in $GL_n(\ZZ)$. Since $\bfmu= E\bfla E^{\tr}$,
the proof of part (b) is complete.

(c) Since $\bfla$ is antisymmetric, the group $\sum_{i,j=1}^n \ZZ \lambda_{ij}$ is
generated by the $\lambda_{ij}$ for $i<j$, so the assumption of rank $n(n-1)/2$
implies that $\{\lambda_{ij} \mid 1\le i<j\le n\}$ is a basis for $\sum_{i,j=1}^n
\ZZ \lambda_{ij}$. As noted in the proof of part (b), $\sum_{i,j=1}^n \ZZ\mu_{ij}=
\sum_{i,j=1}^n \ZZ \lambda_{ij}$, and so this group also has a basis
$\{\mu_{ij} \mid 1\le i<j\le n\}$.

Next, identify $\bfla$ with the linear transformation on $k^n$ given by left
multiplication of $\bfla$ on column vectors. We claim that $\ZZ^n\cap \ker\bfla=
\{\bfz\}$. If $a= (a_1,\dots,a_n)^{\tr} \in  \ZZ^n\cap \ker\bfla$, then
$\lambda_{12}a_2+
\lambda_{13}a_3+ \cdots+ \lambda_{1n}a_n =0$. Since
$\lambda_{12},\dots,\lambda_{1n}$ are $\ZZ$-linearly independent, it follows that
$a_2= a_3= \cdots= a_n= 0$. Then $\lambda_{21}a_1 =0$, which implies $a_1=0$
because $\lambda_{21}= -\lambda_{12} \ne 0$. Thus $a=\bfz$, establishing the
claim. Since $\bfla= (BA)\bfla(BA)^{\tr}$, it follows that $\ZZ^n\cap \ker
(BA)^{\tr}=
\{\bfz\}$. But $(BA)^{\tr}$ is an integer matrix, so we obtain $\det (BA)^{\tr}
\ne 0$, and thus $\det(BA) \ne 0$.

Write $BA= (d_{ij})$, and compare entries in the equation $\bfla=
(BA)\bfla(BA)^{\tr}$:
$$\lambda_{ij}= \sum_{l,m=1}^n d_{il}\lambda_{lm}d_{jm}= \sum_{1\le l<m\le n}
(d_{il}d_{jm}- d_{im}d_{jl}) \lambda_{lm}$$
for all $i$, $j$. Since the $\lambda_{lm}$ for $l<m$ are $\ZZ$-linearly independent,
we find that
$$d_{il}d_{jm}- d_{im}d_{jl}= \delta_{il}\delta_{jm}$$
for $1\le i<j\le n$ and $1\le l<m\le n$. It follows from the Laplace relations that all the
$2\times2$ and larger minors of $BA$ for which the row and column index sets differ must vanish.
In particular, this implies that the adjoint matrix $D= \adj(BA)$ is diagonal. Since $BAD=
\det(BA)I_n$ and $\det(BA)\ne 0$, we conclude that $BA$ must be a diagonal matrix.

The equation $\bfla= (BA)\bfla(BA)^{\tr}$ now reduces to $\lambda_{ij}=
d_{ii}\lambda_{ij}d_{jj}$ for all $i$, $j$, whence $d_{ii}d_{jj}=1$ for all
$i<j$. Since $n\ge 2$ and the $d_{ii}$ are integers, $d_{ii}= \pm1$ for all $i$,
whence $BA\in GL_n(\ZZ)$. Therefore $A\in GL_n(\ZZ)$, proving part (c).
\qed\enddemo

It is tempting to conjecture that the equivalence of Theorem 5.5 holds for arbitrary antisymmetric
$\bfla,\bfmu\in M_n(k)$.

\head Acknowledgements\endhead

We thank J\. Alev for extensive discussions on Poisson
algebras, and T\. Levasseur for helpful correspondence and references concerning
differentials of group actions.

\enddocument